# REGULARIZATION IN KERNEL LEARNING

By Shahar Mendelson[1] and Joseph Neeman

*The Australian National University and Technion,*
*I.I.T. and University of California, Berkeley*

Under mild assumptions on the kernel, we obtain the best known error rates in a regularized learning scenario taking place in the corresponding reproducing kernel Hilbert space (RKHS). The main novelty in the analysis is a proof that one can use a regularization term that grows significantly slower than the standard quadratic growth in the RKHS norm.

**1. Introduction.** Let $F$ be a family of functions from a probability space $(\Omega, \mu)$ to $\mathbb{R}$. A classical problem of learning theory is the following: we set $\nu$ to be an (unknown) probability measure on $\Omega \times \mathbb{R}$ whose marginal distribution on $\Omega$ is $\mu$. Given $n$ independent samples $(X_1, Y_1), \dots, (X_n, Y_n) \in \Omega \times \mathbb{R}$, distributed according to $\nu$, our task is to find a function $\hat{f} \in F$ such that

$$(1.1) \qquad \mathbb{E}(\hat{f}(X_1) - Y_1)^2 - \inf_{f \in F} \mathbb{E}(f(X_1) - Y_1)^2$$

is very small. In other words, we want to approximate the distribution $\nu$ by a function from $F$ as closely as possible. Specifically, we want to find a method of choosing $\hat{f}$ as a function of the sample $(X_i, Y_i)_{i=1}^n$ such that, with high probability, (1.1) is smaller than a function of $n$ that tends to zero as $n$ grows. In this paper, we will consider the case where $\Omega$ is a compact Hausdorff space and $Y_i$ is bounded almost surely.

A widely used approach to solving this problem is to consider a function $\hat{f} \in F$ that minimizes the functional

$$\sum_{i=1}^n (f(X_i) - Y_i)^2$$

Received October 2008; revised April 2009.

[1]Supported in part by Australian Research Council Discovery Grant DP0559465 and by Israel Science Foundation Grant 666/06.

over all $f \in F$. Such a function is called an *empirical minimizer* and its properties have been widely studied (see, e.g., [2, 3, 8, 16, 19] and references therein). It turns out that the complexity and geometry of $F$ play a large part in determining whether (1.1) is small. Roughly speaking, if $F$ is a small family of functions, then (1.1) will be, with high probability, a function of $n$ that decreases polynomially fast.

Of course, there is a disadvantage to having a small family of functions, namely, that $\inf_{f \in F} \mathbb{E}(f(X_1) - Y_1)^2$ becomes larger as $F$ becomes smaller. This trade-off is known as the *bias-variance problem*. The expression (1.1) is known as the *sample error* and $\inf_{f \in F} \mathbb{E}(f(X_1) - Y_1)^2$ is called the *approximation error*.

One major issue that needs to be addressed when using the empirical minimization algorithm is overfitting. Since all of the information that one has is on the behavior of the minimizer on the sample, there is no way to distinguish a "simple" minimizer from a more complicated one. The *regularized learning model* is a method of solving the bias-variance problem while addressing the overfitting problem. We take $F$ to be a very large function class (so that the approximation error is small) and consider a function $\hat{f}$ that minimizes the functional

$$\sum_{i=1}^{n} (f(X_i) - Y_i)^2 + \gamma_n(f),$$

where $\gamma_n(f)$ measures, in some sense, the "complexity" of the function of $f$ and, for a fixed $f$, $\gamma_n(f) \to 0$ as $n \to \infty$. Thus, if two functions have the same empirical behavior, then the algorithm will choose the simpler function of the two.

A common example of the regularized learning problem, and the situation we will be considering in this article, is the case where the class of functions is a reproducing kernel Hilbert space (RKHS), which is defined below and will be denoted throughout this article by $H$. All of the error bounds in this situation (with the exception of one result, discussed later, in the classification setting) were restricted to a regularization term of the form $\gamma_n(f) = \eta_n \|f\|_H^2$, with the goal being to choose $\eta_n$ so that the error is as small as possible. As far as we know, it has not even been conjectured that one could improve the power of $\|f\|_H$ in the regularization process. Doing just that is the main goal of this article.

One can motivate the regularized learning model by looking at it as a collection of empirical minimization problems. Indeed, let $B_H$ be the unit ball of the space $H$ and consider the empirical minimization problem in $rB_H$ for some $r > 0$. As $r$ increases, the approximation error for $rB_H$ decreases and its sample error increases. We could achieve a small total error by choosing the correct value of $r$ and performing empirical minimization in $rB_H$. The role of the regularization term $\gamma_n(f)$ is to force the algorithm to choose the correct value of $r$ for empirical minimization. We will explain later why this motivation can be made rigorous and that the regularization problem may be solved by a solution to a hierarchy of minimization problems.

It should be clear from this motivation that the choice of $\gamma_n$ is critical for the success of the regularized learning model. There has been some significant work done recently on finding explicit formulas for $\gamma_n$ that provide low error rates with high probability. Of particular importance to us, because their results are directly comparable to ours, are the works of Caponnetto and De Vito [7], Smale and Zhou [31] and Wu, Ying and Zhou [36]. We will mention these results in Section 3, in order to compare them to ours. Recent work on regularization parameters for support vector machines includes that of Blanchard, Bousquet and Massart [5]; and Steinwart and Scovel [29]. Although there are important differences between our problem and that of support vector machines, there are certain similarities, and some tools—model selection results and localized complexity parameters, for example—are useful for both subjects.

Our starting point is the realization that analysis based on $L_\infty$-bounds (even if done in a subtle way) is too loose and is the sole source of a quadratic regularization term. Using $L_\infty$-bounds is very tempting in our case because of a convenient fact about reproducing kernel Hilbert spaces: if the kernel is bounded, then there is a constant $c_K$ such that $\|f\|_\infty \le c_K \|f\|_H$. This allows one to bound $\|(f - Y)^2\|_\infty \lesssim \|f\|_H^2$, which can be used to control the "complexity" of the loss class through concentration inequalities (such as Bernstein's for a single function or Talagrand's for a class of functions) that depend on the $L_\infty$-norm of functions. What is more significant is that it allows one to apply contraction inequalities at the cost of a multiplicative factor—the Lipschitz constant of the loss function on its domain—which is, for the squared loss, twice the maximal $L_\infty$-norm of a class member. This approach leads to a regularization term of $\|f\|_H^2$ and it is by avoiding gratuitous use of $L_\infty$-bounds that we improve that term.

It should be noted that in the classification setting (to be more precise, in the example of support vector machines), Blanchard, Bousquet and Massart [5] showed that a regularization term of the form $\eta_n \|f\|_H$ was possible. Their approach unfortunately does not extend to the regression case: they still rely on $L_\infty$-bounds and they obtain a linear regularization term because the loss function in a support vector machine setup is $\ell(x, y) = \max\{0, 1 - xy\}$ and $\|\max\{0, 1 - f(X)Y\}\|_\infty$ is linear in $\|f\|_H$ (instead of quadratic, as is the case for the squared loss). On the other hand, it is conceivable that our technique could be applied to the classification setting, lowering the exponent of $\|f\|_H$ further still.

The starting point of our analysis is the notion of isomorphic coordinate projections, introduced in the context of learning theory in [3]. Suppose that $F$ is a family of functions for which the infimum $\inf_{f \in F} \mathbb{E}(f(X) - Y)^2$ is achieved; call the minimizer $f^*$ and define the excess loss function to be, for any $f \in F$,

$$\mathcal{L}_f^F(X, Y) = (f(X) - Y)^2 - (f^*(X) - Y)^2.$$

When the underlying class is clear from the context, we will omit the superscript $F$. Denote by $P$ the conditional expectation with respect to the sample,

$$P\mathcal{L}_f = \mathbb{E}(\mathcal{L}_f | X_1, Y_1, \ldots, X_n, Y_n)$$

and let $P_n \mathcal{L}_f = \sum_{i=1}^n \mathcal{L}_f(X_i, Y_i)$. One can show (see [3] or Theorem 2.2) that there is some (small) number $\rho_n$ such that, with probability at least $1 - e^{-x}$, every $f \in F$ satisfies

$$(1.2) \qquad \frac{1}{2} P_n \mathcal{L}_f - \rho_n \leq P\mathcal{L}_f \leq 2 P_n \mathcal{L}_f + \rho_n.$$

We will refer to equations like (1.2) as giving "almost isomorphic coordinate projections" because (1.2) tells us that the structures imposed on $F$ by $P$ and $P_n$ are, up to a small additive term, isomorphic. This is a useful approach for bounding the error of the empirical minimizer. Indeed, it is not hard to see that it implies that

$$\mathbb{E}(\hat{f}(X) - Y)^2 - \inf_{f \in F} \mathbb{E}(f(X) - Y)^2 = P\mathcal{L}_{\hat{f}} \leq \rho_n.$$

It turns out that this isomorphic coordinate projection approach applies to regularized learning as well as to empirical minimization. The main result in this direction is due to Bartlett [1] and implies that if every ball $rB_H$ satisfies an almost-isomorphic condition, then it is possible to establish a regularized learning bound. This is an example of a *model selection* result because it proves that the regularized learning procedure somehow selects an appropriate model ($rB_H$ for a good choice of $r$) from a family of models (the set of models $\{rB_H : r \geq 1\}$). Of course, model selection results have been used previously in the study of regularized learning; the use of an almost-isomorphic coordinate projection condition, however, first occurred in [1] and it is crucial here. Some examples of model selection results for problems similar to ours can be found in [5] (Theorem 4.3), [16] and [21].

THEOREM 1.1 [1]. *For each $f \in H$, let $\mathcal{L}_f$ denote the loss of $f$ relative to the ball $\|f\| B_H$:*

$$\mathcal{L}_f(X, Y) = \mathcal{L}_f^{\|f\| B_H}(X, Y) = (f(X) - Y)^2 - (f^*(X) - Y)^2,$$

*where $f^* = \arg\min_{\|g\| \leq \|f\|} \mathbb{E}(g(X) - Y)^2$. Under some conditions on $\gamma_n(\cdot)$, if, for every $f \in H$,*

$$\frac{1}{2} P_n \mathcal{L}_f - \gamma_n(f) \leq P\mathcal{L}_f \leq 2 P_n \mathcal{L}_f + \gamma_n(f),$$

*then the regularized minimizer satisfies*

$$\mathbb{E}(\hat{f}(X) - Y)^2 \leq \inf_{f \in H}((f(X) - Y)^2 + c\gamma_n(c'f)),$$

*where $c$ and $c'$ are absolute constants.*

Thus, if one could establish sharp "isomorphic coordinate projections"-type estimates for every excess loss class $\{\mathcal{L}_f : f \in rB_H\}$, then this would yield regularization bounds.

It is important to emphasize that although at first glance, the problem of obtaining isomorphic bounds for kernel classes has been solved in the past (based on, e.g., estimates from [2, 22]), this is far from being the case. The isomorphic bounds for kernel classes have been studied for the base class $F = B_H$ (i.e., $r = 1$), using an $L_\infty$-based argument that includes contraction inequalities. In contrast, the essential ingredient required for our analysis (and which determines the regularization parameter) is the way in which these bounds scale with the radius $r$. In all of the previous isomorphic results obtained in the context of kernel classes, the way that the bounds depend on $r$ was not important and thus never addressed. And, moreover, the analysis used to obtain those results gives a suboptimal estimate as a function of $r$: an estimate that scales like $r^2$ because of the $L_\infty$-based method. Indeed, one factor of $r$ in this quadratic growth follows from a contraction argument combined with the fact that the maximal $L_\infty$-norm of functions in $rB_H$ is $r$. The second factor of $r$ appears because the "complexity" of the class $rB_H$ grows linearly in $r$.

The consequences of this are clear: since one can identify the regularization term with the way isomorphic coordinate projection estimates for the class $rB_H$ scale with $r$, the regularization term of $\|f\|_H^2$ is an artifact of the $L_\infty$-based method of analysis that leads to a bound that grows like $r^2$.

Let us mention that if the Lipschitz constant of the loss is bounded by an absolute constant, as is the case for support vector machines and the hinge loss, one factor of $r$ can be removed by the $L_\infty$-based method because the Lipschitz constant of the loss is uniformly bounded. Thus, one can use contraction inequalities freely for that problem and obtain a linear regularization term; this is the result in [5].

Our analysis will show that the standard regularization bounds, which grow like $r^2$, where $r = \|f\|_H$, are very pessimistic and may be improved considerably. Moreover, if we set the regularization term as $\eta_n \nu(\|f\|_H)$, we will establish the best known bounds on $\eta_n$ as well (both results will require mild assumptions on the kernel).

There are two reasons for the improved bounds. The first is a method that allows one to bypass the whole $L_\infty$-based mechanism and this is presented in Section 4. We shall present a general bound on the empirical process indexed by the localized excess squared loss class associated with a base class consisting of linear functionals on $\ell_2$ of norm at most $r$. This step will lead to a removal of one factor of $r$ from the $r^2$ term—the one that was due to an $L_\infty$-based method and a contraction argument.

Second, the ability to employ the "isomorphic" approach allows one to use localization techniques. Thus, the effective complexity of the excess loss class is caused only by excess loss functions with a relatively small variance; by virtue of the geometry of $rB_H$, that set of excess loss functions happens to come from a rather small subset of $rB_H$. Recall that, intuitively, the second factor of $r$ comes from the linear growth of the "complexity" of $rB_H$. However, the actual "isomorphic" estimate for $rB_H$ is determined by the complexity of the intersection bodies $xB_2 \cap rB_H$, rather than by that of $rB_H$ [where $B_2$ is the unit ball of $L_2(\mu)$]. It turns out that for a reasonable RKHS, the complexity of such an intersection body grows at a much slower rate as a function of $r$. Indeed, the number of "meaningful directions" in $rB_H$ (when considered as a subset of $L_2$) is small and decreases quickly with $r$. Therefore, the true complexity of $rB_H$ will be sublinear in $r$ because, as $r$ increases, an ever smaller number of directions will actually grow with $r$ and influence the complexity.

Formally, we will show that if the eigenvalues of the integral operator $T_K$ decay like $O(t^{-1/p})$ for some $0 < p < 1$, then one can obtain an isomorphic bound with $\rho_n$ that scales like

$$\max\{\theta^{2/(1+p)}, \theta^{2/p}\}$$

for $\theta \sim r^p n^{-1/2} \log n$. This translates to a regularization term of

$$\max\left\{r^{2p/(1+p)}\left(\frac{\log^2 n}{n}\right)^{1/(1+p)}, \frac{r^2}{n}\right\},$$

where, again, $r = \|f\|_H$.

In this result, one still has a regularization term that grows like $r^2$; nevertheless, this is a considerable improvement on the $L_\infty$-based result. Because it decays faster as a function of the sample size $n$, the $r^2/n$ term seems superfluous because one would expect it to be dominated by the first term. Indeed, in Section 5, we will show that it can be removed: under the same assumption on the decay of the eigenvalues of $T_K$ as above, one may use a regularization term (up to logarithmic term) of

$$\frac{r^{2p/(1+p)}}{n^{1/(1+p)}},$$

which is the best known dependency on $r$ and $n$.

We will end this introduction with the formulation and a short discussion of our main result. To avoid defining them twice, let us mention that the space $\ell_{p,\infty}$ and its norm $\|\cdot\|_{p,\infty}$ are included in Definition 3.3.

ASSUMPTION. Assume that $\|K(x,x)\|_\infty \leq 1$ and that the eigenvalues of the integral operator $T_K$ satisfy $(\lambda_n)_{n=1}^\infty \in \ell_{p,\infty}$ for some $0 < p < 1$. Assume, further, that there is a constant $A$ such that the eigenfunctions $(\varphi_n)_{n\geq 1}$ of $T_K$ satisfy $\sup_n \|\varphi_n\|_\infty \leq A < \infty$.

THEOREM A. *Let $K$ be a continuous, symmetric, positive definite kernel on $\Omega$, a compact Hausdorff space, and set $H$ to be the corresponding reproducing kernel Hilbert space. If $Y$ is bounded almost surely and the assumption above is satisfied, then there exist constants $c_1$, $c_2$ and $c_3$ that depend only on $A$, $p$ and $\|(\lambda_i)\|_{p,\infty}$, a constant $c_Y$ that depends only on $\|Y\|_\infty$ and a constant $N$ depending only on $\|Y\|_\infty$ and $p$ for which the following holds. Let*

$$\tilde{V}(f,u) = c_3(1 + u + c_Y \ln n + \ln \log(\|f\|_H + e)) \left( \frac{(\|f\|_H + 1)^p \log n}{\sqrt{n}} \right)^{2/(1+p)}.$$

*If $n \geq N$ and $c_1 \log \log n \leq u \leq c_2(\log n)^{2/(1-p)}$, then, with probability at least $1 - \exp(-u/2)$, every minimizer $\hat{f}$ of*

$$P_n \ell_f + \kappa_1 \tilde{V}(f,u)$$

*satisfies*

$$P\ell_{\hat{f}} \leq \inf_{f \in H} (P\ell_f + \kappa_2 \tilde{V}(f,u)),$$

*where $\kappa_1$ and $\kappa_2$ are absolute constants and $\ell_f = (f - Y)^2$ is the squared loss function.*

Let us begin our discussion with the assumptions. The assumption that $\|K\|_\infty \leq 1$ is purely cosmetic: any continuous kernel on a compact space is bounded and the assumption only prevents unnecessary constants from appearing. The assumption on the decay of the eigenvalues is essentially a smoothness condition for the kernel; the existence, for example, of a continuous derivative would be enough. We will discuss the eigenvalue assumption in more detail later. For now, let us just say that it has been used before [7] in discussing the way in which smoothness of the kernel affects the learning rates.

The assumption that the eigenfunctions $\varphi_n$ are uniformly bounded is more serious. It has been made before—in [35], for example, in which it was mistakenly claimed that such an assumption holds for all Mercer kernels. Zhou, [37], however, argues against this assumption and provides an example of a $C^\infty$ kernel without uniformly bounded eigenfunctions. Let us remark, therefore, that we do not need the full strength of this assumption. Indeed, as the proof will reveal, it is enough to have some $0 < \epsilon < 1/2$ such that $\sup_n \lambda_n^\epsilon \|\varphi_n\|_\infty$ is bounded. The theorem then remains true if we assume that $(\lambda_n^{1-2\epsilon}) \in \ell_{p,\infty}$ instead of $(\lambda_n) \in \ell_{p,\infty}$. Note that $\sup_n \sqrt{\lambda_n} \|\varphi_n\|_\infty < \infty$; for our assumption to hold, we need to be able to take a power of $\lambda_n$ that is strictly smaller than $1/2$. This is a considerably weaker assumption than that of uniformly bounded eigenfunctions. For instance, the example given in [37] of a $C^\infty$ kernel without uniformly bounded eigenfunctions satisfies our weaker condition for any $\epsilon > 0$: the eigenvalues decrease exponentially faster than the $L_\infty$-norms of the eigenfunctions.

For an example of a kernel satisfying our assumption, let $k$ be an even function of period 1 and set $K(x, y) = k(x - y)$. If $\mu$ is the Lebesgue measure on $[0, 1]$, then it is easily seen, via a cosine expansion of $k$, that the eigenfunctions of $K$ are sine and cosine functions and hence bounded uniformly. The periodic Gaussian kernel is an example of such a kernel.

As a final remark on the assumptions, let us point out that one can trivially construct examples of kernels that satisfy them: just take $\varphi_n$ to be a suitably smooth orthonormal basis of $L_2(\mu)$ and choose $\lambda_n$ to be a sequence that decreases sufficiently rapidly. Then $K(x, y) = \sum_{n=1}^\infty \lambda_n \varphi_n(x) \varphi_n(y)$ satisfies our assumptions.

Regarding the theorem, there are some aspects of practical interest that we do not address. First, there are constants in the theorem that we have made no attempt to compute. Furthermore, these constants depend on quantities that may not be known (e.g., $\|Y\|_\infty$). One might hope, however, to use applied statistical techniques—cross-validation, for example—to find plausible values for the constants. In that case, one should note that the constant $c_Y$ in the definition of $\tilde{V}$ can be moved to the front of the definition without changing the validity of the theorem (see Remark 2.6); that way, the applied statistician has only one unknown constant to worry about.

We conclude this introduction with a brief discussion of the error rate of Theorem A; more detailed discussions follow at the ends of Sections 3 and 5. The formulation of Theorem A is attractive because it shows that we find the almost-minimizer (in some sense) regardless of how well our hypothesis class approximates the regression function $\mathbb{E}(Y|X)$. To be concrete, however, we can make an assumption about how the approximation error behaves and derive explicit error bounds as a function of $n$. The assumption made in [9] (and elsewhere) is that there exists some $0 < \sigma \le 1/2$ such that $\mathbb{E}(Y|X)$ is in the range of $T_K^\sigma$ on $L_2(\mu)$. For $\sigma = 1/2$, this implies that $\mathbb{E}(Y|X) \in H$; for smaller $\sigma$, it somehow says that $\mathbb{E}(Y|X)$ can be approximated reasonably well by elements in $H$. Under this assumption, we obtain an error rate of (ignoring logarithmic factors and the confidence term, $u$) $n^{-2\sigma/(p+2\sigma)}$. As stated above, a detailed discussion follows in later sections; for now, we will just mention that the above rate is significantly faster than the rate of $n^{-\sigma/2}$ that was obtained in [31].

Regarding the optimality of this error rate, we have very little to say. Minimax lower bounds on the error rate are given in [7], but only when the regression function $\mathbb{E}(Y|X)$ belongs to $H$ (and their proof does not easily extend to the more general case considered here). In a very specific case (when $\sigma = 1/2$ and one cannot take $\sigma > 1/2$), our rates match those in [7]. We can claim, therefore, that our results are optimal in a very specific sense; in the more interesting region $0 < \sigma < 1/2$, however, we cannot make any such claim.

**2. Preliminaries.** We begin with a word about notation. We will denote absolute constants (i.e., fixed, positive numbers) by $c, c_1, \ldots$, etc. Their values may change from line to line. Absolute constants whose values will remain unchanged are denoted by $\kappa_1, \kappa_2, \ldots$. By $c(a)$, we mean that the constant $c$ depends only on the parameter $a$. We write $a \sim b$ if there exist absolute constants $c_1$ and $c_2$ such that $c_1 a \leq b \leq c_2 a$, and $a \sim_p b$ if the equivalence constants depend on the parameter $p$.

Arguably the most important tool in modern empirical processes theory is Talagrand's concentration inequality for an empirical process indexed by a class of uniformly bounded functions [18, 33]. The version of this concentration result which we shall use here is due to Massart [20].

THEOREM 2.1. *There exists an absolute constant $C$ for which the following holds. Let $F$ be a class of functions defined on $(\Omega, \mu)$ such that for every $f \in F$, $\|f\|_\infty \leq b$ and $\mathbb{E}f = 0$. Let $X_1, \ldots, X_n$ be independent random variables distributed according to $\mu$ and set $\sigma^2 = n \sup_{f \in F} \mathbb{E}f^2$. Define*

$$Z = \sup_{f \in F} \sum_{i=1}^n f(X_i) \quad and \quad \bar{Z} = \sup_{f \in F} \left| \sum_{i=1}^n f(X_i) \right|.$$

*Then, for every $x > 0$ and every $\rho > 0$,*

$$\Pr(\{Z \geq (1+\rho)\mathbb{E}Z + \sigma\sqrt{Cx} + C(1+\rho^{-1})bx\}) \leq e^{-x},$$

$$\Pr(\{Z \leq (1-\rho)\mathbb{E}Z - \sigma\sqrt{Cx} - C(1+\rho^{-1})bx\}) \leq e^{-x}$$

*and the same inequalities hold for $\bar{Z}$.*

Throughout this article, we denote by $\ell(x, y) = (x-y)^2$ the squared loss function. When $f$ is a function $\Omega \to \mathbb{R}$ and $Y$ is some target random variable, we define $\ell_f = \ell_f(X, Y) = (f(X) - Y)^2$. If $F$ is a class of functions, let $\mathcal{L}_f^F = \mathcal{L}_f(X, Y) = (f(X) - Y)^2 - (f^*(X) - Y)^2$, where $f^* = \arg\min_{f \in F} \mathbb{E}\ell_f$ (we will usually drop the superscript $F$). Of course, we assume that this minimizer exists and is unique, which is the case, for example, if $F$ is compact (in $L_2$) and convex. $\mathcal{L}_F$ denotes the class of functions $\{\mathcal{L}_f^F : f \in F\}$.

For a class of functions $F$ on a probability space $(\Omega, \mu)$, we set

$$\|P_n - P\|_F = \sup_{f \in F} \left| \frac{1}{n} \sum_{i=1}^n f(X_i) - \mathbb{E}f \right|,$$

where $(X_i)_{i=1}$ are independent, distributed according to $\mu$.

For any $x \geq 0$, define the localized excess loss class

$$\mathcal{L}_x = \{\mathcal{L}_f : \mathbb{E}\mathcal{L}_f \leq x\}$$

and set

$$V = \text{star}(\mathcal{L}_F, 0) = \{\theta \mathcal{L}_f : 0 \leq \theta \leq 1, f \in F\},$$

$$V_x = \{\theta \mathcal{L}_f : 0 \leq \theta \leq 1, \mathbb{E}(\theta \mathcal{L}_f) \leq x\} = \{h \in \text{star}(\mathcal{L}_F, 0) : \mathbb{E}h \leq x\}$$

[where, for a set $T$, $\text{star}(T, 0) = \{\theta t : 0 \leq \theta \leq 1, t \in T\}$ is the star-shaped hull of $T$ and $0$].

The following "isomorphic" result is similar in nature to the one proved in [3]. The bound from Theorem 2.2 normally leads to an estimate on the error of the empirical minimizer, but in [4] and here, it will serve a different purpose. This isomorphic result will enable us to control the solution of the regularized learning problem in the context of kernel learning.

THEOREM 2.2. *There exists an absolute constant $c$ for which the following holds. Let $\mathcal{L}_F$ be a squared loss class associated with a convex class $F$ and a random variable $Y$. If $b = \max\{\sup_{f \in F} \|f\|_\infty, \|Y\|_\infty\}$ and $x > 0$ satisfies*

$$\mathbb{E}\|P_n - P\|_{V_x} \leq x/8,$$

*then, with probability $1 - \exp(-u)$, for every $f \in F$,*

(2.1)    $\dfrac{1}{2} P_n \mathcal{L}_f - \dfrac{x}{2} - c(1 + b^2)\dfrac{u}{n} \leq P\mathcal{L}_f \leq 2P_n \mathcal{L}_f + \dfrac{x}{2} + c(1 + b^2)\dfrac{u}{n}.$

PROOF. By Talagrand's inequality, there exists an absolute constant $C$ such that, for every $\alpha > 0$, with probability at least $1 - e^{-u}$,

$$\|P_n - P\|_{V_\alpha} \leq 2\mathbb{E}\|P_n - P\|_{V_\alpha} + \left(\frac{Cu}{n}\right)^{1/2} \sup_{g \in V_\alpha} \sqrt{\text{Var } g} + \frac{Cbu}{n}.$$

It is standard to verify (see, e.g., [19]), that there exists an absolute constant $C$ such that, for a convex class $F$, every $\mathcal{L}_f \in \mathcal{L}_F$ satisfies $\mathbb{E}\mathcal{L}_f^2 \leq Cb^2 \mathbb{E}\mathcal{L}_f$. Thus, every $g \in V_\alpha$ satisfies $\text{Var } g \leq Cb^2 \alpha$. Fix $x$ satisfying $\mathbb{E}\|P_n - P\|_{V_x} \leq x/8$ and set

$$\alpha = \max\left\{x, 25C\frac{(1 + b^2)u}{n}\right\}.$$

Note that, because $V$ is star shaped, $\alpha \geq x$ implies that $V_\alpha \subset \frac{\alpha}{x}V_x$ and so $\mathbb{E}\|P_n - P\|_{V_\alpha} \leq \frac{\alpha}{x}\mathbb{E}\|P_n - P\|_{V_x} \leq \alpha/8$. Therefore, with probability at least

$1 - e^{-u}$,

$$\|P_n - P\|_{V_\alpha} \le \frac{\alpha}{4} + \left( C \frac{b^2 \alpha u}{n} \right)^{1/2} + \frac{Cbu}{n}$$

(2.2)
$$\le \frac{\alpha}{4} + \frac{\alpha}{5} + \frac{\alpha}{25}$$

$$\le \frac{\alpha}{2}.$$

Consider the event in which (2.2) holds. Fix some $\mathcal{L}_f \in \mathcal{L}_F$. If $P\mathcal{L}_f \le \alpha$, then $\mathcal{L}_f \in V_\alpha$ and so

$$P_n \mathcal{L}_f - \frac{\alpha}{2} \le P\mathcal{L}_f \le P_n \mathcal{L}_f + \frac{\alpha}{2},$$

and (2.1) holds. If, on the other hand, $P\mathcal{L}_f = \beta > \alpha$, then let $g = \frac{\alpha}{\beta}\mathcal{L}_f$ and note that $g \in V_\alpha$. Thus, by (2.2),

$$\frac{1}{2}Pg = Pg - \frac{\alpha}{2} \le P_n g \le Pg + \frac{\alpha}{2} \le 2Pg.$$

Since $\mathcal{L}_f$ is a constant multiple of $g$, we have

$$\frac{1}{2}P\mathcal{L}_f \le P_n \mathcal{L}_f \le 2P\mathcal{L}_f$$

and so (2.1) holds once again.

To conclude, (2.2) implies that (2.1) holds for all $\mathcal{L}_f \in \mathcal{L}_F$. Thus, (2.1) holds with probability at least $1 - e^{-u}$. $\square$

REMARK 2.3. The claim of Theorem 2.2 holds under milder assumptions. Note that the assumption that $F$ is convex is there to ensure that $P\ell_f$ attains a unique minimum in $F$ and that the excess loss class satisfies a Bernstein-type condition: that for every $f \in F$, $\mathbb{E}\mathcal{L}_f^2 \le C\mathbb{E}\mathcal{L}_f$. One can show that if $F$ is convex, then, for any function $f \in F$, $\mathbb{E}\mathcal{L}_f^2 \le c\|f\|_\infty^2 \mathbb{E}\mathcal{L}_f$. Hence, if $F$ is convex and $G$ is a subset of $F$ that contains the minimizer in $F$ of $P\ell_f$, then the analog of Theorem 2.2 will be true for $\{\mathcal{L}_g : g \in G\}$.

The first part of our analysis will be to show that this isomorphic information can be used to derive estimates in regularized learning.

2.1. *From isomorphic information to regularized learning.* The regularized learning model provides a method for learning in a very large class of functions without suffering a large statistical error. As we mentioned in the Introduction, obtaining an "isomorphic" result for a hierarchy of classes can lead to estimates in the regularized learning model. This approach was introduced in [1] and was formulated in the way we will use here in [4]. Since this last article has not yet appeared, we present a proof of the result we need in the Appendix.

Let $F$ be a class of functions and suppose that there is a collection of subsets $\{F_r : r \geq 1\}$ with the following properties:

1. $\{F_r : r \geq 1\}$ is monotone (i.e., whenever $r \leq s$, $F_r \subseteq F_s$);
2. for every $r \geq 1$, there exists a unique element $f_r^* \in F_r$ such that $P\ell_{f_r^*} = \inf_{f \in F_r} P\ell_f$;
3. the map $r \rightarrow P\ell_{f_r^*}$ is continuous;
4. for every $r_0 \geq 1$, $\bigcap_{r > r_0} F_r = F_{r_0}$;
5. $\bigcup_{r \geq 1} F_r = F$.

DEFINITION 2.4. Given a class of functions $F$, we say that $\{F_r ; r \geq 1\}$ is an *ordered, parameterized hierarchy of $F$* if the above conditions 1–5 are satisfied. Define, for $f \in F$,

$$r(f) = \inf\{r \geq 1 ; f \in F_r\}.$$

Note that, from the semicontinuity property of an ordered, parameterized hierarchy (property 4), it follows that $f \in F_{r(f)}$ for all $f \in F$.

From the second property of an ordered, parameterized hierarchy, we can, for $r \geq 1$ and $f \in F_r$, define $\mathcal{L}_{r,f} = (f - Y)^2 - (f_r^* - Y)^2$. That is, $\mathcal{L}_{r,f}$ is the excess loss function with respect to the class $F_r$.

THEOREM 2.5. *There exist absolute constants $\kappa_1$ and $\kappa_2$ such that the following holds. Suppose that $\{F_r ; r \geq 1\}$ is an ordered, parameterized hierarchy and that $\rho_n(r, u) : [1, \infty) \times (0, \infty) \rightarrow (0, \infty)$ is a continuous function (possibly depending on the sample) that is increasing in both $r$ and $u$. Suppose, also, that for every $r \geq 1$ and every $u > 0$, with probability at least $1 - \exp(-u)$,*

$$\tfrac{1}{2} P_n \mathcal{L}_{r,f} - \rho_n(r, u) \leq P \mathcal{L}_{r,f} \leq 2 P_n \mathcal{L}_{r,f} + \rho_n(r, u)$$

*for all $f \in F_r$.*

*Then, for every $u > 0$, with probability at least $1 - \exp(-u)$, any function $\hat{f} \in F$ that minimizes the functional*

$$P_n \ell_f + \kappa_1 \rho_n(2r(f), \theta(r(f), u))$$

*also satisfies*

$$P\ell_{\hat{f}} \leq \inf_{f \in F}(P\ell_f + \kappa_2 \rho_n(2r(f), \theta(r(f), u))),$$

*where*

$$\theta(r, x) = x + \ln \frac{\pi^2}{6} + 2 \ln\left(1 + \frac{P\ell_{f_1^*}}{\rho_n(1, x + \log(\pi^2/6))} + \log r\right).$$

REMARK 2.6. In fact, the proof of Theorem 2.5 reveals something slightly stronger: if $\tilde{\rho}_n(r, u)$ is a continuous, increasing function in both variables such that

$$\tilde{\rho}_n(r, u) \geq \rho_n(2r, \theta(r, u))$$

for every $r$, $u$ and $n$, then every function $\hat{f}$ that minimizes the functional

$$P_n\ell_f + \kappa_1\tilde{\rho}_n(r, u)$$

satisfies

$$P\ell_{\hat{f}} \leq \inf_{f \in F}(P\ell_f + \kappa_2\tilde{\rho}_n(r, u)).$$

In other words, we can always regularize with a larger regularization term; we will obtain a correspondingly larger error bound. We will use this fact later.

The conclusion of Theorem 2.5 can be reformulated in a way that makes the traditional distinction between the approximation and sample errors more explicit. We begin by defining an approximation error term by

$$\mathcal{A}(r) = \inf_{f \in F_r} P\ell_f.$$

Then $\mathcal{A}(r) - \inf_{f \in F} P\ell_f$ tends to zero as $r \to \infty$ and the rate of this convergence measures how well the ordered, parameterized hierarchy approximates $Y$. Smale and Zhou [30] study this approximation error in a variety of contexts, including the case in which we are interested: when $F_r$ is the ball of radius $r - 1$ in a reproducing kernel Hilbert space.

COROLLARY 2.7. *Under the assumptions of Theorem 2.5, with probability at least $1 - \exp(-u)$,*

$$P\ell_{\hat{f}} \leq \inf_{r \geq 1}(\mathcal{A}(r) + \kappa_2\rho_n(2r, \theta(r, u))).$$

PROOF. Let $u > 0$, fix $\varepsilon > 0$ and choose an $s \geq 1$ such that

$$\mathcal{A}(s) + \kappa_2\rho_n(2s, \theta(s, u)) \leq \inf_{r \geq 1}(\mathcal{A}(r) + \kappa_2\rho_n(2r, \theta(r, u))) + \frac{\varepsilon}{2}.$$

Consider $g \in F_s$ such that $P\mathcal{L}_g \leq A(s) + \varepsilon/2$. Since $\rho_n$ is increasing in both of its arguments, we have

$$P\mathcal{L}_g + \kappa_2\rho_n(2r(g), \theta(r(g), u)) \leq \inf_{r \geq 1}(\mathcal{A}(r) + \kappa_2\rho_n(2r, \theta(r, u))) + \varepsilon.$$

However, we can find such a function $g$ for every $\varepsilon > 0$. Therefore,

$$\inf_{f \in F}(P\mathcal{L}_f + \kappa_2\rho_n(2r(f), \theta(r(f), u))) \leq \inf_{r \geq 1}(\mathcal{A}(r) + \kappa_2\rho_n(2r, \theta(r, u)))$$

and the conclusion follows from Theorem 2.5. □

**3. Regularization in kernel classes.** The case that we will be interested in is when $F_r$ is a multiple of the unit ball of an RKHS. For more details on properties of an RKHS that are relevant in the context of learning theory, we refer the reader to, for example, [8].

Let $\Omega$ be a compact Hausdorff space, consider $K:\Omega \times \Omega \to \mathbb{R}$, a positive definite, continuous function and, without loss of generality, assume that $\|K\|_\infty \leq 1$. Let $T_K$ be the corresponding integral operator, $T_K:L_2(\mu) \to L_2(\mu)$, defined by

$$(T_K f)(x) = \int_\Omega K(x,y)f(y)\,d\mu(y).$$

By Mercer's theorem [8], there is an orthonormal basis of eigenfunctions $(\varphi_i)_{i=1}^\infty$ of $T_K$, corresponding to the eigenvalues $(\lambda_i)_{i=1}^\infty$ arranged in a non-increasing order, such that

$$K(x,y) = \sum_{i=1}^\infty \lambda_i \varphi_i(x) \varphi_i(y),$$

where the convergence is uniform and absolute on the support of $\mu \times \mu$ [10] (and hence there is also convergence in $L_2$).

The RKHS, which will be denoted throughout by $H$, can be identified with linear functionals in $\ell_2$. Indeed, consider the function $\Phi:\Omega \to \ell_2$ defined by $\Phi(x) = (\sqrt{\lambda_i}\varphi_i(x))_{i=1}^\infty$. For every $t \in \ell_2$, define the corresponding element of $H$ by $f_t(x) = \langle \Phi(x), t \rangle$; we define the RKHS $H$ to be the image of $\ell_2$ under the map $t \mapsto f_t$ with the induced inner product $\langle f_t, f_s \rangle_H = \langle t, s \rangle$. This definition of $H$ is phrased differently from those given in [8, 10], but it is easily checked that the resulting Hilbert space of functions is the same. Hence, to study properties of a subset of $H$, it is enough to study the corresponding set of linear functionals, as a set $T \subset \ell_2$ uniquely determines $F_T = \{f_t : t \in T\}$. Here, we will mostly be concerned with $T = rB_2$, corresponding to $F = rB_H$, where $B_H$ is the unit ball of the RKHS and $B_2$ is the unit ball of $\ell_2$. In this case, the measure endowed on $\ell_2$ is given by $\Phi(Z)$, where $Z$ is distributed in $\Omega$ according to $\mu$.

3.1. *Classes of linear functionals: The $L_\infty$ approach.* Our first approach to the problem of regularized learning in an RKHS will lead to a regularization term of $\|f\|_H^2$. As stated in the Introduction, this is over-regularization, which is an artifact of the analysis of the learning problem. It stems from the way that the $L_\infty$-bound on functions in $\mathcal{L}_F$ is used and the fact that the only way to bound $\|\mathcal{L}_f\|_{L_\infty}$ is by $\|\mathcal{L}_f\|_{L_\infty} \leq c\|f\|_H^2$. In this section, we will use this (loose) approach, but still obtain better error estimates than those previously known—although still using a regularization term of $\|f\|_H^2$. We will obtain considerably better results in the following sections.

The idea we will use is to obtain an isomorphic result for the hierarchy $F_r = rB_H$ (in our $\ell_2$ representation, $F_r$ corresponds to $rB_2$). We then use Corollary 2.7 for the function $\rho_n$ given by the isomorphic analysis.

In our presentation, we will study the following, more general, situation. Let $T \subset \ell_2$ be a compact, convex, symmetric set and consider a random vector $\xi$ on $\ell_2$ [distributed, recall, according to $\Phi(Z)$]. Denote by $f_t = \langle t, \cdot \rangle$ the linear functional defined by $t$ and put

$$D = \{t : \mathbb{E}f_t(X)^2 \le 1\} = \{t : \mathbb{E}\langle t, \xi \rangle^2 \le 1\}.$$

Thus, $D$ is the image of the $L_2$ unit ball in the parameter space $\ell_2$.

Our first, $L_\infty$-based, approach to the problem of learning in an RKHS relies on the following bound, which was implicit in [22] (to be precise, Theorem 3.1 follows from the proof of Theorem 3.3 in [22] if one keeps track explicitly of the constants $C_b$ and $C_b'$).

THEOREM 3.1. *There exist constants $c$ and $c'$ depending only on $\|Y\|_\infty$ for which the following holds. Let $V_{r,x} = \{\alpha\mathcal{L}_f; 0 \le \alpha \le 1, f \in rB_H, \mathbb{E}\mathcal{L}_f \le x\}$. Then for every $r \ge 1$ and every $x > 0$,*

$$\mathbb{E}\|P - P_n\|_{V_{r,x}} \le cr\mathbb{E}\sup_{\{t \in rB_2 \cap \sqrt{x}D\}}\left|\frac{1}{n}\sum_{i=1}^n g_i f_t(X_i)\right|,$$

*where the $g_i$ are independent standard Gaussian variables. In the case where $r = 1$, we have*

$$\mathbb{E}\sup_{\{t \in B_2 \cap \sqrt{x}D\}}\left|\sum_{i=1}^n g_i f_t(X_i)\right| \le c'\left(\frac{1}{n}\sum_{i=1}^\infty \min\{x, \lambda_i\}\right)^{1/2}.$$

The proof of the first part of Theorem 3.1 uses a comparison theorem, relating the Gaussian process $t \to \sum_{i=1}^n g_i\mathcal{L}_{f_t}(X_i, Y_i)$, conditioned on $(X_i, Y_i)_{i=1}^n$, to the conditioned Gaussian process $t \to \sum_{i=1}^n g_i f_t(X_i)$. This is done using an $L_\infty$-bound since

$$\sum_{i=1}^n (\mathcal{L}_{f_t} - \mathcal{L}_{f_s})^2(X_i, Y_i) = \sum_{i=1}^n (f_t - f_s)^2(X_i) \cdot ((f_t + f_s)(X_i) - 2Y_i)^2$$

$$\le 4(r + \|Y\|_\infty)^2 \sum_{i=1}^n (f_t - f_s)^2(X_i),$$

which will turn out to be the main source of the quadratic regularization term $\|f\|_H^2$.

From Theorem 3.1, one obtains the following.

COROLLARY 3.2. *There exists a constant $\tilde{c}$, depending only on $\|Y\|_\infty$, such that if $z > 0$ satisfies*

$$z \geq \tilde{c}\left(\frac{1}{n}\sum_{i=1}^{\infty}\min\{z, \lambda_i\}\right)^{1/2},$$

*then, for all $r \geq 1$,*

$$\frac{x}{8} \geq \mathbb{E}\|P - P_n\|_{V_{r,x}},$$

*where $x = r^2 z$.*

PROOF. Define

$$\psi_r(x) = r\mathbb{E}\sup_{\{t \in rB_2 \cap \sqrt{x}D\}}\left|\sum_{i=1}^{n}g_i f_t(X_i)\right|.$$

By the second part of Theorem 3.1, we can choose $\tilde{c}$ such that $\psi_r(x) \leq \frac{x}{8c}$ (where $c$ is the constant from Theorem 3.1). Furthermore, it is easily checked that $\psi_r(x) = r^2\psi_1(xr^{-2})$ for any $x$ and $r$. That is, $\psi_r(r^2x) = r^2\psi_1(x) \leq \frac{r^2x}{8c}$. The claim now follows from the first part of Theorem 3.1. □

With this corollary and Theorem 2.2, we can obtain an isomorphic condition on the unit ball of an RKHS using information on the decay of the eigenvalues. For the sake of concreteness, we will make the following assumption on this rate of decay; this assumption will allow us to compute an error bound explicitly.

DEFINITION 3.3. For $0 < p < 1$ and a nonincreasing, nonnegative sequence $(\lambda_i)_{i=1}^{\infty}$, define

$$\|(\lambda_i)\|_{p,\infty} = \sup_{i \geq 1} i^{1/p}\lambda_i.$$

Hence, for any $x > 0$,

(3.1)                    $|\{\lambda_i \geq x\}| \leq \|(\lambda_i)\|_{p,\infty}x^{-p}.$

If $\|(\lambda_i)\|_{p,\infty} < \infty$, we will say that $(\lambda_i) \in \ell_{p,\infty}$.

ASSUMPTION 3.1. Let $K$ be a kernel on a compact probability space $(\Omega \times \Omega, \mu \times \mu)$ where $\mu$ is a Borel measure and $\Omega \subset \mathbb{R}^d$. Assume that $\|K(x,x)\|_\infty \leq 1$ and that the eigenvalues of the integral operator $T_K$ satisfy $(\lambda_n)_{n=1}^{\infty} \in \ell_{p,\infty}$ for some $0 < p < 1$.

Since $\int K(x,x)\,d\mu(x) = \sum_{i=1}^{\infty} \lambda_i$, we have $(\lambda_i) \in \ell_{1,\infty}$ when $K(x,x) \in L_1(\mu)$. The stronger Assumption 3.1 is satisfied under some smoothness condition on the kernel. Suppose, for example, that the kernel $K$ belongs to some Besov space $B_{2,\infty}^{\alpha}$ [in particular, this is the case if $\alpha \in \mathbb{N}$ and $K \in \mathcal{C}^{\alpha}(\Omega \times \Omega)$]. If $\Omega \subset \mathbb{R}^d$ is locally the graph of a Lipschitz function and $\mu$ is a Borel (probability) measure on $\Omega$, then, by Theorems 4.1 and 4.7 of [6] (see also [17]), the sequence $(\lambda_i)$ belongs to $\ell_{p,\infty}$ for

$$p = \frac{1}{\alpha/d + 1/2}.$$

A similar assumption on the decay of the eigenvalues was made in [7]. The $L_\infty$ assumption on $K(x,x)$ is only to simplify the presentation and any uniform bound instead of 1 would do.

The assumption on the rate of decay of the eigenvalues allows us to obtain the following bound.

LEMMA 3.4. *For $0 < p < 1$, there is a constant $c_p$ depending only on $p$ such that for all $x > 0$ and all $r > 0$,*

$$\sum_{i=1}^{\infty} \min\{x, r^2 \lambda_i\} \le c_p \|(\lambda_i)\|_{p,\infty} x^{1-p} r^{2p}.$$

PROOF. It suffices to prove the lemma for $r = 1$ and the result will follow for all $r$ by homogeneity. Set $N_x = |\{\lambda_i \ge x\}|$ and observe that for all $x > 0$,

$$\sum_{i=1}^{\infty} \min\{x, \lambda_i\} = xN_x + \sum_{i=N_x+1}^{\infty} \lambda_i \le \|(\lambda_i)\|_{p,\infty} x^{1-p} + \sum_{i=N_x+1}^{\infty} \lambda_i.$$

The first term is in the required form. Let us deal with the second term:

$$\sum_{i=N_x+1}^{\infty} \lambda_i \le \|(\lambda_i)\|_{p,\infty} \sum_{i=N_x+1}^{\infty} i^{-1/p}$$

$$\le c_p \|(\lambda_i)\|_{p,\infty} N_x^{1-1/p}$$

$$\le c_p \|(\lambda_i)\|_{p,\infty} x^{p-1}$$

as required. □

With the preceding bound on $\frac{1}{n}\sum_i \min\{x, r^2\lambda_i\}$, we can rewrite Corollary 3.2 in a nicer form that is specialized to our application; recall that $V_{r,x}$ is the localization at level $x$ of the star-shaped hull of the shifted loss class of $rB_H : V_{r,x} = \{\alpha \mathcal{L}_f : 0 \le \alpha \le 1, f \in rB_H, \mathbb{E}\mathcal{L}_f \le x\}$.

COROLLARY 3.5.  *Let $K$ be a kernel that satisfies Assumption 3.1 for some $0 < p < 1$. There exists a constant $c_p$ depending only on $p$ such that if $z = c_p (\frac{\|(\lambda_i)\|_{p,\infty}}{n})^{1/(1+p)}$, then, for all $r > 1$,*

$$\frac{x}{8} \geq \mathbb{E}\|P - P_n\|_{V_{r,x}},$$

*where $x = r^2 z$.*

Having controlled the quantity, $\mathbb{E}\|P - P_n\|_{V_{r,x}}$, that can give us the "isomorphic coordinate projection" result we wanted, we are almost in a position to prove our first result; it only remains to show that we can apply the model selection result.

LEMMA 3.6.  *Let $H$ be the reproducing kernel Hilbert space associated with a continuous, symmetric, positive definite kernel $K$. Set $F = H$ and define, for every $r \geq 1$, $F_r = (r-1)B_H$, where $B_H$ is the closed unit ball of $H$. Then $\{F_r; r \geq 1\}$ is an ordered, parameterized hierarchy and $r(f) = \|f\| + 1$.*

PROOF.  The first, fourth and fifth properties of an ordered, parameterized hierarchy are immediate. The second property follows from the fact that $B_H$ is convex and compact with respect to the $L_2$-norm (because it is an ellipsoid whose principal lengths decrease to zero). For the third property, fix $1 \leq q < r < s$ and let $\beta = \frac{q-1}{r-1}$, $\alpha = \frac{r-1}{s-1}$. Note that $\alpha f_s^* \in F_r$ and $\beta f_r^* \in F_q$. Thus,

$$0 \leq P\ell_{f_r^*} - P\ell_{f_s^*} \leq P\ell_{\alpha f_s^*} - P\ell_{f_s^*} = (\alpha^2 - 1)P(f_s^*)^2 + 2(1-\alpha)Pf_s^* Y.$$

As $s \to r$, the right-hand side tends to zero (because the candidates for $f_s^*$ are uniformly bounded in $L_2$) and so $r \to P\ell_{f_r^*}$ is upper semicontinuous (the same argument works for $r = 1$). In the other direction,

$$0 \leq P\ell_{f_q^*} - P\ell_{f_r^*} \leq P\ell_{\beta f_r^*} - P\ell_{f_r^*} \leq (\beta^2 - 1)P(f_r^*)^2 + 2(1-\beta)Pf_r^* Y$$

and the right-hand side tends to zero for the same reason as before.  □

Combining Theorem 2.2 with Corollaries 3.5 and 2.7, we obtain the following error bound for regularized learning in an RKHS.

THEOREM 3.7.  *There exist absolute constants $\kappa_1$ and $\kappa_2$, constants $c_Y$ and $c_Y'$ depending only on $\|Y\|_\infty$ and a constant $c_p$ depending only on $p$ such that the following holds. Let $K$ be a kernel satisfying Assumption 3.1 and define*

$$\rho_n(r, u) = c_p r^2 \left(\frac{\|(\lambda_i)\|_{p,\infty}}{n}\right)^{1/(1+p)} + c_Y(1 + r^2)\frac{u}{n}.$$

*Then, for every $u > 0$, with probability at least $1 - \exp(-u)$, any function $\hat{f} \in F$ that minimizes the functional*

$$P_n \ell_f + \kappa_1 \tilde{\rho}_n(r(f), u)$$

*also satisfies*

$$P\ell_{\hat{f}} \leq \inf_{r \geq 1}(\mathcal{A}(r) + \kappa_2 \tilde{\rho}_n(r, u)),$$

*where*

$$\tilde{\rho}_n(r, u) = \rho_n\left(2r, u + \ln\frac{\pi^2}{6} + 2\ln(1 + c'_Y n + \log r)\right).$$

*In particular,*

$$P\ell_{\hat{f}} \leq \inf_{r \geq 1}\left(\mathcal{A}(r) + c\left(\frac{r^2}{n^{1/(1+p)}} + \frac{1 + r^2}{n}(u + \log n + \log\log(r + e))\right)\right),$$

*where $c = c(p, \|Y\|_\infty, \|(\lambda_i)\|_{p,\infty})$.*

PROOF.    By Theorem 2.2 and Corollary 3.5, the function $\rho_n(r, x)$ satisfies the condition of Theorem 2.5 [where we set $c_Y = c(1 + \|Y\|_\infty^2)$]. We can then apply Corollary 2.7 to obtain the result. Since $0 \in F_r$ for any $r > 0$, we have $P\ell_{f_1^*} \leq P\ell_0 = \|Y\|_{L_2(\mu)}^2$ and $\rho_n(1, u + \ln(\pi^2/6)) \geq c''_Y/n$ so that

$$\frac{P\ell_{f_1^*}}{\rho_n(1, x + \ln(\pi^2/6))} \leq c'_Y n,$$

to which we apply Remark 2.6.    □

Let us compare the estimate on the regularization term and the resulting error rate that follows from this theorem to previously obtained bounds on regularized learning in an RKHS. Since all of the results we consider have exponentially good confidence, we will simplify this comparison by ignoring the confidence term and focusing on the decay of the error bound as the sample size increases. In order to facilitate our comparison further, we will make an assumption that allows us to control the approximation error $\mathcal{A}(r)$.

ASSUMPTION 3.2.    Suppose that there exists $0 < \sigma < 1$ such that $T_K^{-\sigma}\mathbb{E}(Y|X)$ belongs to $L_2$.

Recall that $T_K$ is the integral operator that defines our RKHS $H$. Note that for $\sigma = 0$, the assumption is trivial and for $\sigma \geq \frac{1}{2}$, the assumption states that $\mathbb{E}(Y|X) \in H$ [and so $\mathcal{A}(r) = 0$ for large enough $r$]. For $0 < \sigma < \frac{1}{2}$, the assumption tells us the degree to which $\mathbb{E}(Y|X)$ can be approximated by functions in $H$. Indeed, a result of Smale and Zhou [30] (see also [10]) allows us to bound the approximation error in terms of $\sigma$, as follows.

THEOREM 3.8 [30]. *If Assumption 3.2 holds for $0 < \sigma < \frac{1}{2}$, then*

$$\mathcal{A}(r-1) - \inf_{f \in H} P\ell_f \leq \left(\frac{1}{r}\right)^{4\sigma/(1-2\sigma)} \|T_K^{-\sigma} \mathbb{E}(Y|X)\|_2^{2/(1-2\sigma)}.$$

Our main points of comparison are the rates from [31], Corollary 5:

$$(3.2) \qquad P\ell_{\hat{f}} - \inf_{f \in H} P\ell_f \lesssim \begin{cases} \left(\dfrac{1}{n}\right)^{\sigma/(1+2\sigma)}, & \text{if } \sigma \geq \frac{1}{2}, \\[2mm] \left(\dfrac{1}{n}\right)^{\sigma/2}, & \text{if } \sigma < \frac{1}{2}. \end{cases}$$

Suppose, first, that Assumption 3.2 holds for $\sigma \geq \frac{1}{2}$. As we already mentioned, this implies that the approximation error is eventually zero and so our result gives an error rate like

$$P\ell_{\hat{f}} - \inf_{f \in H} P\ell_f \lesssim \left(\frac{1}{n}\right)^{1/(1+p)},$$

which is an improvement over (3.2), even if $p = 1$. In fact, [7] shows that this rate is optimal in some sense. Interestingly, [7] also shows that one can get even better rates for $\sigma > \frac{1}{2}$, that is, when the regression function not only belongs to the hypothesis class, but also satisfies some extra smoothness properties.

For $\sigma < \frac{1}{2}$, set $k = 4\sigma/(1-2\sigma)$. We can then choose $r = n^{1/((1+p)(2+k))}$ and our error bound becomes

$$P\ell_{\hat{f}} - \inf_{f \in H} P\ell_f \lesssim \mathcal{A}(n^{1/((1+p)(2+k))}) + n^{2/((1+p)(2+k))} \left(\frac{1}{n}\right)^{1/(p+1)}$$

$$(3.3) \qquad\qquad \lesssim \left(\frac{1}{n}\right)^{2\sigma/(1+p)}.$$

Once again, this improves on (3.2), even when $p = 1$.

The situation $p < 1$ is more interesting because the kernels used in learning theory often have some smoothness properties. If $K \in C^\infty$, for example, then we can choose $p$ arbitrarily small and recover the following result of Wu, Ying and Zhou [36]:

$$P\ell_{\hat{f}} - \inf_{f \in H} P\ell_f \lesssim \left(\frac{1}{n}\right)^{2\sigma - \epsilon}$$

for any $\epsilon > 0$. We will see, however, that the techniques of the next section will improve on this for $\sigma < \frac{1}{2}$.

**4. Toward a smaller regularization parameter.** The bound (3.3) would be substantially improved if we could remove the $r^2$ term and replace it by a smaller power of $r$—which is the main source of novelty in this article. As mentioned before, the most significant source for this improvement comes from bypassing $L_\infty$-based bounds. In recent years, there has been considerable progress made on bounding various empirical processes that are indexed by sets that are either not bounded or very weakly bounded in $L_\infty$. Most of these results were motivated by questions in asymptotic geometric analysis, most notably, sampling from an isotropic, log-concave measure (e.g., [15, 25, 28]) and the approximate reconstruction problem [14, 23]. The fact that such an approach is called for here seems strange because we are dealing with a learning problem relative to a class of uniformly bounded functions, so it would seem that there is no reason to employ techniques designed to handle an unbounded situation. Even more so, because in a standard learning analysis, the way the error bounds depend on the $L_\infty$-diameter of the class is usually of no real importance. In contrast, here, the way the isomorphic results scale with the $L_\infty$-bound is extremely important because one is trying to obtain a result for the entire hierarchy and the $L_\infty$-diameter of $F_r$ is directly linked to the hierarchy parameter $r$. Thus, the standard, and very loose, approach which is commonly used in a single class situation can cause real damage in our case because the regularization term will be strongly influenced by the way that the $L_\infty$-diameter enters into the bounds.

To see where one can improve upon the standard $L_\infty$ analysis (in a very "hand-waving" way), let us return to the localized Gaussian process indexed by $\{t : \mathbb{E}\mathcal{L}_{f_t} \leq x\} \cap rB_2$, conditioned on the data $(X_i, Y_i)$, that is,

$$t \to \sum_{i=1}^{n} g_i \mathcal{L}_{f_t}(X_i, Y_i) = \sum_{i=1}^{n} g_i \langle t - t^*, X_i \rangle (\langle t + t^*, X_i \rangle - 2Y_i),$$

where $f_{t^*}$ minimizes the loss in $rB_2$. For every $t$, the variance of each conditioned Gaussian variable satisfies

$$\sigma^2 \left( \sum_{i=1}^{n} g_i \mathcal{L}_{f_t}(X_i, Y_i) \right) = \sum_{i=1}^{n} \langle t - t^*, X_i \rangle^2 (\langle t + t^*, X_i \rangle - 2Y_i)^2.$$

Consider some $t$ for which $\mathbb{E}\mathcal{L}_{f_t} \leq x$. One can show that in this case, $\|t - t^*\| \leq \sqrt{x}$ (see Lemma 4.1 below). Now, if one has a very strong concentration phenomenon and if $D = B_2$, then

$$\left( \left\langle \frac{t + t^*}{2}, X_i \right\rangle - Y_i \right)^2 = \left( \left\langle \frac{t - t^*}{2}, X_i \right\rangle + (\langle t^*, X_i \rangle - Y_i) \right)^2$$
$$\approx_c x + \mathbb{E}\ell_{f_{t^*}}.$$

Since the expected loss of the best in the class only decreases with $r$, this term is of the order of $x$, rather than a factor that grows quadratically in $r$, which is the estimate that results from the $L_\infty$ approach. This at least hints at the fact that the $L_\infty$ approach is likely to lead to very loose estimates.

Despite the fact that the above paragraph is totally unjustified as stated and very optimistic, it turns out that this scenario is very close to the actual situation (although the proof requires a rather delicate analysis).

4.1. *Further preliminaries.* For technical reasons, we will make an additional assumption on the eigenfunctions of the kernel. We should emphasize that it is possible that this assumption may not be necessary to obtain the improved regularization term, although we were not able to remove it here and it has a crucial role in our analysis.

ASSUMPTION 4.1. Let $K$ be a kernel on a compact probability space $(\Omega \times \Omega, \mu \times \mu)$ with $\Omega \subset \mathbb{R}^d$. Assume that there is a constant $A$ such that the eigenfunctions of $K$ satisfy $\sup_n \|\varphi_n\|_\infty \leq A < \infty$.

Let us recall from the Introduction that we still obtain a result if we assume instead that there exists $\epsilon > 0$ with $\sup_n \lambda_n^\epsilon \|\varphi_n\|_\infty \leq A < \infty$. As discussed in the Introduction, our results hold with this weaker assumption if we modify Assumption 3.1 so that $(\lambda_n^{1-2\epsilon})_{n=1}^\infty \in \ell_{p,\infty}$.

Recall that the feature map $\Phi$ defines an isometry from an RKHS into $\ell_2$. Let $T \subset \Phi(H)$ be a centrally symmetric, convex, compact subset of $\ell_2$. The first step in our analysis is to relate the localized sets $\mathcal{L}_x$ (corresponding to the class $\{f_t : t \in T\}$) to subsets of $T$. Since this fact appeared implicitly in several places (see, e.g., [24], Corollary 3.4) and in more general situations, for example, loss functions that are uniformly convex rather than the squared loss, we omit its proof.

LEMMA 4.1. *Let $t^* = \arg\min_{t \in T} \mathbb{E}\ell_{f_t}$. For every $x > 0$,*

$$\{t - t^* : t \in T, \mathcal{L}_{f_t} \in \mathcal{L}_x\} \subset 2\sqrt{x} D \cap 2T.$$

Lemma 4.1 shows that it is sufficient to consider the complexity of the sets $\sqrt{x} D \cap T$. The complexity parameters we shall use come from a generic chaining argument (defined below) and thus a significant part of our analysis will be based on covering numbers.

DEFINITION 4.2. Let $A, B \subset \ell_2$. Denote by $N(A, B)$ the smallest number of translates of $B$ needed to cover $A$. If $\varepsilon B$ is a ball of radius $\varepsilon$ with respect to some norm, then $N(A, \varepsilon B)$ is the minimal cardinality of an $\varepsilon$-cover of $A$ with respect to that norm. If $(A, d)$ is a metric space (rather than a normed one), we denote the cardinality of a minimal $\varepsilon$-cover of $A$ by $N(A, \varepsilon, d)$.

The generic chaining mechanism (see [34] for the most recent survey on this topic) is used to relate probabilistic properties of a random process indexed by a metric space to the metric structure of the underlying space. This mechanism originated in the study of Gaussian processes $t \to X_t$, where it was proven that $\mathbb{E} \sup_{t \in T} X_t$ is equivalent to a metric invariant of $(T, d)$ for $d(s, t) = (\mathbb{E}|X_s - X_t|^2)^{1/2}$. This so-called *majorizing measures theorem* (in which the upper bound of the equivalence was proven by Fernique [12] and the lower by Talagrand [32]) was later developed into a more general theory with many interesting applications [34]. The metric invariant that is at the heart of this theory is the $\gamma_2$ functional, which we define as follows.

Let $(T, d)$ be a metric space. An *admissible sequence* of $T$ is a collection of subsets of $T$, $\{T_s : s \geq 0\}$, such that for every $s \geq 1$, $|T_s| = 2^{2^s}$ and $|T_0| = 1$.

DEFINITION 4.3. For a metric space $(T, d)$, define

$$\gamma_2(T, d) = \inf \sup_{t \in T} \sum_{s=0}^{\infty} 2^{s/2} d(t, T_s),$$

where the infimum is taken with respect to all admissible sequences of $T$ and $d(t, T) = \inf_{u \in T} d(t, u)$.

DEFINITION 4.4. A random process $t \to X_t$ indexed by a metric space $(T, d)$ is sub-Gaussian relative to $d$ if, for every $s, t \in T$ and every $u \geq 1$,

$$\Pr(|X_s - X_t| \geq u \, d(s, t)) \leq 2 \exp\left(-\frac{u^2}{2}\right).$$

The generic chaining mechanism can be used to show that if $\{X_t : t \in (T, d)\}$ is sub-Gaussian, then there is an absolute constant $c$ such that for every $t_0 \in T$,

$$\mathbb{E} \sup_{t \in T} |X_t - X_{t_0}| \leq c \gamma_2(T, d)$$

and similar bounds hold with high probability.

Note that one choice for sets $T_s$ that constitute a potential (yet, usually suboptimal) admissible sequence are $\varepsilon_s$-covers of $T$, where each $\varepsilon_s$ is selected in a way that ensures that $N(T, \varepsilon_s, d) \leq 2^{2^s}$. An easy computation [34] then shows that

$$(4.1) \qquad \gamma_2(T, d) \leq c \int_0^{\operatorname{diam}(T, d)} \sqrt{\log N(T, \varepsilon, d)} \, d\varepsilon,$$

where $c$ is an absolute constant. This is a generalization of Dudley's entropy integral (see, e.g., [11, 34]), used in the study of Gaussian processes. As will be explained later, this integral bound can be improved under certain assumptions on the geometry of $T$ if $d$ is endowed with a norm.

The metric $d$ we will focus on here is a random one and depends on the sample $X_1, \ldots, X_n \subset \ell_2$. For every $X_1, \ldots, X_n$, set

$$d_{\infty,n}(f, g) = \max_{1 \leq i \leq n} |f(X_i) - g(X_i)|.$$

Recall that our function class $H$ is isometric to $\ell_2$ under the map $t \mapsto f_t$. Thus, $d_{\infty,n}$ defines a random norm on a projection of $\ell_2$ which is given, with some abuse of notation, by

$$d_{\infty,n}(s, t) = \max_{1 \leq i \leq n} |\langle X_i, s - t \rangle|.$$

Next, let $U_n(T) = (\mathbb{E}\gamma_2^2(T, d_{\infty,n}))^{1/2}$ and, for every $x > 0$, set

$$\phi_n(x) = \frac{U_n(K_x)}{\sqrt{n}} \cdot \max\left( \sqrt{x}, \sqrt{\mathbb{E}\mathcal{L}_{t^*}}, \frac{U_n(K_x)}{\sqrt{n}} \right),$$

where $K_x = T \cap \sqrt{x}D \subset \ell_2$ and $t^*$ is the parameter in $T$ for which $\inf_{t \in T} \mathbb{E}\mathcal{L}_{f_t}$ is attained.

Recall that

$$\mathcal{L}_x = \{\mathcal{L}_f : \mathbb{E}\mathcal{L}_f \leq x\}$$

and that

$$V_x = \{\theta\mathcal{L}_f : 0 \leq \theta \leq 1, \mathbb{E}(\theta\mathcal{L}_f) \leq x\} = \{h \in \text{star}(\mathcal{L}_F, 0) : \mathbb{E}h \leq x\}.$$

From Theorem 2.2, it is clear that in order to obtain a useful "isomorphic" result, one has to bound $\mathbb{E}\|P_n - P\|_{V_x}$ as a function of $x$; this is done in the following theorem. Since it is a modification of a result that was proven in [4], we will only present an outline of its proof.

THEOREM 4.5. *There exists an absolute constant $c$ for which the following holds. If $T \subset \ell_2$ and $H = \{f_t : t \in T\}$, then, for every $x > 0$,*

$$\mathbb{E}\|P_n - P\|_{V_x} \leq c\sum_{i=0}^{\infty} 2^{-i}\phi_n(2^{i+1}x).$$

The proof of Theorem 4.5 relies on the following "peeling" lemma, which shows that one can control $\mathbb{E}\|P_n - P\|$ on the star-shaped hull of a class of functions if one can control $\mathbb{E}\|P_n - P\|$ on "shells" of the original class.

LEMMA 4.6. *For every $x > 0$,*

$$\mathbb{E}\|P_n - P\|_{V_x} \leq 2\sum_{i=0}^{\infty} 2^{-i}\mathbb{E}\|P_n - P\|_{\mathcal{L}_{2^{i+1}x}}.$$

PROOF. Note that for every $x > 0$,

$$
\begin{aligned}
W_x &= \{\theta \mathcal{L}_f : 0 \le \theta \le 1, \mathbb{E}(\theta \mathcal{L}_f) \le x, \mathbb{E}\mathcal{L}_f \ge x\} \\
&= \left\{ \frac{t\mathcal{L}_f}{\mathbb{E}\mathcal{L}_f} : \mathbb{E}\mathcal{L}_f \ge x, 0 \le t \le x \right\} \\
&= \bigcup_{i=0}^{\infty} \left\{ \frac{t\mathcal{L}_f}{\mathbb{E}\mathcal{L}_f} : 2^i x \le \mathbb{E}\mathcal{L}_f \le 2^{i+1}x, 0 \le t \le x \right\} \equiv \bigcup_{i=0}^{\infty} W_{i,x}.
\end{aligned}
$$

If $t\mathcal{L}_f / \mathbb{E}\mathcal{L}_f \in W_{i,x}$, then $t/\mathbb{E}\mathcal{L}_f \le 2^{-i}$ and $\mathcal{L}_f \in \mathcal{L}_{2^{i+1}x}$. Thus, $\|P_n - P\|_{W_{i,x}} \le 2^{-i}\|P_n - P\|_{\mathcal{L}_{2^{i+1}x}}$.

Finally, let $W_{0,x} = \operatorname{star}(\mathcal{L}_x, 0)$. Note that $\|P_n - P\|_{W_{0,x}} \le \|P_n - P\|_{\mathcal{L}_x}$ and that $V_x \subset W_0 \cup W_{0,x}$, from which our claim follows. $\quad\square$

OUTLINE OF THE PROOF OF THEOREM 4.5. Fix $x > 0$. First, one can verify that the Bernoulli process indexed by $\mathcal{L}_x$, given by $t \to \sum_{i=1}^{n} \varepsilon_i \mathcal{L}_{f_t}(X_i, Y_i)$ conditioned on $(X_i, Y_i)_{i=1}^{n}$ is sub-Gaussian with respect to the metric

$$
d(f_{t_1}, f_{t_2}) = d_{\infty,n}(f_{t_1}, f_{t_2}) \left( \sup_{v \in \sqrt{x}D \cap T} \sum_{i=1}^{n} \langle X_i, v \rangle^2 + \sum_{i=1}^{n} \mathcal{L}_{t^*}(X_i, Y_i) \right)^{1/2}.
$$

This follows from Hoeffding's inequality [which says that the process is sub-Gaussian with respect to $d(\mathcal{L}_{f_t}, \mathcal{L}_{g_t})$] and a computation to show that $d(\mathcal{L}_{f_t}, \mathcal{L}_{g_t})$ is smaller than the above quantity. Hence, if we set $K = \sqrt{x}D \cap T$, then, by the Giné–Zinn symmetrization method [13], followed by a generic chaining argument, we have

$$
\mathbb{E}\|P_n - P\|_{\mathcal{L}_x} \le \frac{c_1}{n} \mathbb{E} \left( \gamma_2(K, d_{\infty,n}) \left( \sup_{t \in K} \sum_{i=1}^{n} \langle t, X_i \rangle^2 + \sum_{i=1}^{n} \mathcal{L}_{t^*}(X_i, Y_i) \right)^{1/2} \right).
$$

Moreover, one can show (see, e.g., [14]) that if $H$ is a class of functions, then

$$
\mathbb{E} \sup_{h \in H} \left| \sum_{i=1}^{n} h^2(X_i) - \mathbb{E}h^2 \right| \le c_2 \max\{\sqrt{n}\sigma_H U_n(H), U_n^2(H)\},
$$

where $\sigma_H^2 = \sup_{h \in H} \mathbb{E}h^2$. In particular, for $H = \{\langle t, \cdot \rangle : t \in K\}$,

$$
\mathbb{E} \sup_{t \in K} \sum_{i=1}^{n} \langle t, X_i \rangle^2 \le nx + c_2 \max\{\sqrt{nx} U_n(K), U_n^2(K)\},
$$

because $\mathbb{E}\langle t, \cdot \rangle^2 \le x$. Now, a straightforward computation shows that

$$
\mathbb{E}\|P_n - P\|_{\mathcal{L}_x} \le \phi_n(x).
$$

To conclude the proof, note that by Lemma 4.6, it is possible to estimate $\mathbb{E}\|P_n - P\|_{V_x}$ using $\mathbb{E}\|P_n - P\|_{\mathcal{L}_{2^i x}}$. $\quad\square$

Observe that the sets $T$ we will be interested in are $rB_2$ since they are the images of $rB_H$ in $\ell_2$. The rest of this section will be devoted to finding a bound on $\phi_n(x)$ for these sets $T$.

4.2. *Controlling $\phi_n$ for $T = rB_2$.* It is clear that $\phi_n$ is determined by the structure of the sets $K_{x,r} = \sqrt{x}D \cap 2rB_2 \subset \ell_2$. To study the metric properties of these sets, we first have to identify $D$.

Consider the random variable $Z$ on $\Omega$ distributed according to $\mu$ and let $X = \Phi(Z) = \sum_{i=1}^{\infty} \sqrt{\lambda_i}\varphi_i(Z)e_i \in \ell_2$ be the random feature map. Clearly,

$$D = \{t \in \ell_2 : \mathbb{E}\langle t, X\rangle^2 \leq 1\} = \{t \in \ell_2 : \mathbb{E}\langle t, \Phi(Z)\rangle^2 \leq 1\}.$$

Since $(\varphi_i)_{i=1}^{\infty}$ is an orthonormal system in $L_2(\mu)$, we have

$$\mathbb{E}\langle t, \Phi(Z)\rangle^2 = \mathbb{E}\sum_{i,j} t_i t_j \sqrt{\lambda_i \lambda_j}\varphi_i(Z)\varphi_j(Z) = \sum_{i=1}^{\infty}\lambda_i t_i^2.$$

Hence, $D$ is an ellipsoid in $\ell_2$ with the standard basis $(e_i)_{i=1}^{\infty}$ as principal directions, and lengths $1/\sqrt{\lambda_i}$.

It is straightforward to verify that for every $x, r > 0$, there is an ellipsoid $\mathcal{E}_{x,r}$ such that $K_{x,r} = 2rB_2 \cap \sqrt{x}D$ satisfies $\frac{1}{2}\mathcal{E}_{x,r} \subset K_{x,r} \subset \mathcal{E}_{x,r}$. The principal directions of $K_{x,r}$ and $\mathcal{E}_{x,r}$ coincide and the principal lengths of $\mathcal{E}_{x,r}$ are

$$c\min\left\{\sqrt{\frac{x}{\lambda_i}}, r\right\},$$

where $c$ is an absolute constant.

The structure of the ellipsoids $\mathcal{E}_{x,r}$ indicates that it should be possible to obtain a sublinear dependency on the radius $r$ and the fact that we were not able to do so in Section 3.1 is an artifact of the suboptimal analysis that was used there. The sublinearity occurs because for $\alpha > 1$, $\mathcal{E}_{x,\alpha r}$ is much smaller than $\alpha\mathcal{E}_{x,r}$; since it is an intersection body, it only grows in some directions and the number of directions in which it grows decreases quickly with $r$.

Now that we have identified the intersection body, we are ready to estimate

$$U_n = (\mathbb{E}\gamma_2^2(\mathcal{E}_{x,r}, d_{\infty,n}))^{1/2}.$$

THEOREM 4.7. *There exists an absolute constant $c$ for which the following holds. Suppose $\sup_n \|\varphi_n\|_\infty \le A$ and set*

$$Q(x, r) = A\left(\sum_{i=1}^\infty \min\{x, r^2\lambda_i\}\right)^{1/2}.$$

*Then*

$$(\mathbb{E}\gamma_2^2(\mathcal{E}_{x,r}, d_{\infty,n}))^{1/2} \le cQ(x, r)\log n.$$

Before proving the theorem, we need two additional facts. The first is an improved "Dudley entropy integral" bound, due to Talagrand.

THEOREM 4.8 [34]. *There exists an absolute constant $c$ for which the following holds. If $\mathcal{E} \subset \mathbb{R}^m$ is an ellipsoid and $B$ is the unit ball of some norm $\|\cdot\|$ on $\mathbb{R}^m$, then*

$$\gamma_2(\mathcal{E}, \|\cdot\|) \le c\left(\int_0^\infty \varepsilon \log N(\mathcal{E}, \varepsilon B)\, d\varepsilon\right)^{1/2}.$$

Another standard fact we need is the dual Sudakov inequality [26].

LEMMA 4.9. *There exists an absolute constant $c$ for which the following holds. Let $B_E$ be the unit ball of some norm on $\mathbb{R}^m$ and let $B_2^m$ be the Euclidean ball on $\mathbb{R}^m$. Then, for every $\varepsilon > 0$,*

$$\log N(B_2^m, \varepsilon B_E) \le c\left(\frac{\mathbb{E}\|G\|_E}{\varepsilon}\right)^2,$$

*where $G = (g_1, \ldots, g_m)$ is a standard Gaussian vector on $\mathbb{R}^m$.*

PROOF OF THEOREM 4.7. Fix $X_1, \ldots, X_n$ and note that in order to bound $\gamma_2(\mathcal{E}_{x,r}, d_{\infty,n})$, it suffices to consider the projection of the (infinite-dimensional) ellipsoid $\mathcal{E}_{x,r}$ onto the subspace spanned by $X_1, \ldots, X_n$. Hence, one can apply Lemma 4.9. Set $\|v\|_E = \max_{1 \le i \le n} |\langle v, X_i\rangle|$ and let $B_E$ be the unit ball $\{v \in \ell_2 : \|v\|_E \le 1\}$. Consider the ellipsoid $\mathcal{E}_{x,r} \subset \ell_2$ with principal directions $(e_i)_{i=1}^\infty$ and lengths $\theta_i = c_1 \min\{\sqrt{x/\lambda_i}, r\}$. Let $T$ be the operator $Te_i = \theta_i e_i$ so that $TB_2 = \mathcal{E}_{x,r}$. For every $\varepsilon > 0$,

$$N(TB_2, \varepsilon B_E) = N(B_2, \varepsilon T^{-1}B_E)$$

and $v \in \varepsilon T^{-1}B_E$ if and only if $\max_{1 \le i \le n} |\langle v, T^*X_i\rangle| = \max_{1 \le i \le n} |\langle v, TX_i\rangle| \le \varepsilon$. Hence, if we set $W_i = TX_i$ and $\|v\|_{\bar{E}} = \max_{1 \le i \le n} |\langle v, W_i\rangle|$ (with the corresponding unit ball $B_{\bar{E}} = \{v : \|v\|_{\bar{E}} \le 1\}$), then

$$N(TB_2, \varepsilon B_E) = N(B_2, \varepsilon B_{\bar{E}}) = N(B_2^n, \varepsilon B_{\bar{E}}),$$

where, here, by $B_2^n$, we mean the unit ball in the subspace of $\ell_2$ spanned by $(W_i)_{i=1}^n$.

Let $G$ be a standard Gaussian vector on $\mathbb{R}^n$. Then, by Slepian's lemma [11, 27],

$$\mathbb{E}\|G\|_{\bar E} = \mathbb{E} \max_{1 \le i \le n} |\langle G, TX_i \rangle| \le c_2 \sqrt{\log n} \max_{1 \le i \le n} \|TX_i\|_2.$$

Since $T$ is a diagonal operator and $X_j = \sum_{i=1}^\infty \sqrt{\lambda_i}\varphi_i(Z_j)e_i$, we have

$$\|TX_j\|_2^2 = \sum_{i=1}^\infty \theta_i^2 \lambda_i \varphi_i^2(Z_j) \le A^2 \sum_{i=1}^\infty \theta_i^2 \lambda_i = A^2 \sum_{i=1}^\infty \min\{x, r^2\lambda_i\}.$$

Hence, setting

$$Q = Q(x,r) = A\left(\sum_{i=1}^\infty \min\{x, r^2\lambda_i\}\right)^{1/2},$$

it is evident that

(4.2) $$\mathbb{E}\|G\|_{\bar E} \le c_2 \sqrt{\log n}\, Q$$

and by Lemma 4.9, for every $\varepsilon > 0$,

$$\log N(B_2^n, \varepsilon B_{\bar E}) \le c_3 \frac{Q^2 \log n}{\varepsilon^2}.$$

In particular, the diameter of $B_2^n$ with respect to the norm $\|\cdot\|_{\bar E}$ is at most $cQ\sqrt{\log n}$, and we denote this diameter by $D_2$.

This estimate for the covering numbers will be used for "large" scales of $\varepsilon$. For smaller scales, we need a different argument. Applying a volumetric estimate (see, e.g., [27]) for every norm $\|\cdot\|_X$ on $\mathbb{R}^n$ and every $\varepsilon > 0$, we have $N(B_X, \varepsilon B_X) \le (5/\varepsilon)^n$. Thus, for every $0 < \varepsilon < \delta$,

$$\log N(B_2^n, \varepsilon B_{\bar E}) \le \log N(B_2^n, \delta B_{\bar E}) + \log N(\delta B_{\bar E}, \varepsilon B_{\bar E})$$
$$\le c_3 \frac{Q^2 \log n}{\delta^2} + n \log\left(\frac{\delta}{\varepsilon}\right).$$

If we take $\delta^2 = c_3 Q^2 \frac{\log n}{n}$, then it follows that for $\varepsilon \le c_4 Q \sqrt{\log n / n} = \varepsilon_0$,

$$\log N(B_2^n, \varepsilon B_{\bar E}) \le n \log(\varepsilon_0/\varepsilon).$$

Now, by Theorem 4.8, for every $X_1, \ldots, X_n$,

$$\gamma_2^2(\mathcal{E}_{x,r}, d_{\infty,n}) \le c_5 \int_0^\infty \varepsilon \log N(TB_2, \varepsilon B_E)\, d\varepsilon = c_5 \int_0^\infty \varepsilon \log N(B_2^n, \varepsilon B_{\bar E})\, d\varepsilon$$
$$\le c_6 \int_0^{\varepsilon_0} n\varepsilon \log\left(\frac{\varepsilon_0}{\varepsilon}\right) d\varepsilon + c_6 \int_{\varepsilon_0}^{D_2} \frac{Q^2 \log n}{\varepsilon}\, d\varepsilon.$$

Using the change of variables $\eta = \varepsilon/\varepsilon_0$, the first integral is bounded by $c_6 n \varepsilon_0^2 \int_0^1 \eta \log(\eta^{-1})\, d\eta = c_7 Q^2 \log n$. Noting that $\varepsilon_0 = c_8 D_2 n^{-1/2}$, the second integral is just

$$c_7 Q^2 \log n (\log D_2 - \log \varepsilon_0) = c_7 Q^2 \log n (\tfrac{1}{2} \log n - \log c_8) \le c_9 Q^2 \log^2 n. \quad \square$$

We will now bound $\phi_n(x)$ using a parameter that describes the decay of the eigenvalues $(\lambda_i)$. By Assumption 3.1, the sequence of eigenvalues has a bounded weak $\ell_p$-norm for some $0 < p < 1$, implying that for all $x > 0$,

$$(4.3) \qquad |\{\lambda_i \ge x\}| \le \|(\lambda_i)\|_{p,\infty} x^{-p}.$$

Set $\tilde{Q}^2(x,r) = c_p A^2 x^{1-p} r^{2p} \|(\lambda_i)\|_{p,\infty}$ and define the function $\tilde{U}_n(x,r)$ by

$$\tilde{U}_n(x,r) = c'_p \tilde{Q}(x,r) \log n,$$

where $c'_p$ is an appropriate constant that depends only on $p$. Then, by Lemma 3.4, $U_n(\mathcal{E}_{x,r}) \le \tilde{U}_n(x,r)$ and setting

$$\tilde{\phi}_n(x,r) = \frac{\tilde{U}_n(x,r)}{\sqrt{n}} \cdot \max\left(\sqrt{x}, \sqrt{\mathbb{E}\mathcal{L}_{t^*}}, \frac{\tilde{U}_n(x,r)}{\sqrt{n}}\right),$$

it follows that for $T = r B_2$, we have $\phi_n(x) \le \tilde{\phi}_n(x,r)$.

LEMMA 4.10. *Suppose that $K$ satisfies Assumptions 3.1 and 4.1. There then exists a constant $c_p$, depending only on $p$, for which the following holds. Let $T_r = r B_2$ and set $V_r$ to be the star-shaped hull of $\{\mathcal{L}_f : f \in T_r\}$. If $V_{r,x} = \{\mathcal{L}_f \in V_r : \mathbb{E}\mathcal{L}_f \le x\}$, then*

$$\mathbb{E}\|P_n - P\|_{V_{r,x}} \le c_p \tilde{\phi}_n(x,r).$$

PROOF. In view of Theorem 4.5, it is enough to show that the sum

$$\sum_{i=0}^{\infty} 2^{-i} \tilde{\phi}_n(2^{i+1} x, r)$$

is dominated by a multiple of the first term in the sum.

For any $\alpha \ge 1$ and any $x > 0$, it is evident from the definition of $\tilde{U}_n$ that

$$\tilde{U}_n(\alpha x, r) \le \alpha^{1/2 - p/2} \tilde{U}_n(x,r);$$

therefore, one can verify that $\tilde{\phi}_n(\alpha x, r) \le \alpha^{1-p/2} \tilde{\phi}_n(x,r)$. In particular,

$$\sum_{i=0}^{\infty} 2^{-i} \tilde{\phi}_n(2^{i+1} x, r) \le 2^{1-p/2} \sum_{i=0}^{\infty} 2^{-ip/2} \tilde{\phi}_n(x,r) \le c_p \tilde{\phi}_n(x,r). \qquad \square$$

Let us pause and explain why this analysis indeed yields a far better result than the $L_\infty$ approach. We will show later that the dominant factor in $\mathbb{E}\|P_n - P\|_{V_{r,x}}$ is $\tilde{U}_n/\sqrt{n}$, which is, up to a logarithmic term and appropriate constants,

$$A\left(\frac{1}{n}\sum_{i=1}^{\infty}\min\{x, r^2\lambda_i\}\right)^{1/2} = (*).$$

In comparison, the $L_\infty$ approach leads to a bound of the order of

$$r\left(\frac{1}{n}\sum_{i=1}^{\infty}\min\{x, r^2\lambda_i\}\right)^{1/2} = (**)$$

on $\mathbb{E}\|P_n - P\|_{V_{r,x}}$—which is considerably larger as $r$ tends to infinity.

If $x$ is a "fixed point" of $(**)$ (as required in the "isomorphic" result of Theorem 2.2), then

$$\left(\frac{1}{n}\sum_{i=1}^{\infty}\min\left\{\frac{x}{r^2}, \lambda_i\right\}\right)^{1/2} = c\frac{x}{r^2}$$

and thus $x$ scales quadratically in $r$. On the other hand, the fixed point of $(*)$ satisfies

$$rA\left(\frac{1}{n}\sum_{i=1}^{\infty}\min\left\{\frac{x}{r^2}, \lambda_i\right\}\right)^{1/2} = cx.$$

Hence, if $(\lambda_i)$ decays quickly, then the fixed point will scale like a smaller power of $r$—in the worst case, linearly in $r$.

The estimate on the fixed point in the alternative approach we presented in this section is the following.

THEOREM 4.11.    *There exists a constant $c_{p,Y}$ depending only on $p$ and $\|Y\|_{L_2}$ such that the following holds. If Assumptions 3.1 and 4.1 are satisfied, then for every $r > 1$, if*

$$\Theta = \frac{A\|(\lambda_i)\|_{p,\infty}^{1/2}r^p\log n}{\sqrt{n}}$$

*and*

$$x \geq c_{p,Y}\max\{\Theta^{2/(1+p)}, \Theta^{2/p}\},$$

*then one has*

$$\mathbb{E}\|P_n - P\|_{V_{x,r}} \leq x/8.$$

PROOF. Fix $r > 1$. From the definition of $\tilde{\phi}_n$, it suffices to find $x$ for which $\tilde{U}_n(x, r)/\sqrt{n} \leq c_Y \min\{x, \sqrt{x}\}$, where $c_Y \leq c_1 \min\{1, (\mathbb{E}\mathcal{L}_{t^*})^{-1/2}\}$, for a suitable absolute constant $c_1$. Note that since $t = 0$ is a potential minimizer, $c_Y \leq c_1(1 + (\mathbb{E}Y^2)^{1/2})$.

The definition of $\Theta$ ensures that $\tilde{U}_n(x, r)/\sqrt{n} = c'_p x^{1/2 - p/2} \Theta$. To have $\tilde{U}_n(x, r)/\sqrt{n} \leq cx$, therefore, it is enough to have $x \geq (c_{p,Y}\Theta)^{2/(1+p)}$. Similarly, to have $\tilde{U}_n(x, r)/\sqrt{n} \leq cx^{1/2}$, it is enough that $cx \geq (c_{p,Y}\Theta)^{2/p}$. $\quad\square$

COROLLARY 4.12. *Suppose that Assumptions 3.1 and 4.1 hold. There then exists a constant $c_{p,Y,A,\lambda}$ depending on $p$, $\|Y\|_\infty$, $\|(\lambda_i)\|_{p,\infty}^{1/2}$ and $A$ such that the function $\rho_n$ defined by*

$$\rho_n(r, u) = c_{p,Y,A,\lambda}(1 + u) \max\left\{\frac{r^{2p/(1+p)} \log^{2/(1+p)}}{n^{1/(1+p)}}, \frac{r^2}{n}\right\}$$

*satisfies the hypothesis of Theorem 2.5.*

*In particular, for every $u > 0$, with probability at least $1 - \exp(-u)$, any function $\hat{f} \in F$ that minimizes the functional*

$$P_n \ell_f + \kappa_1 \tilde{\rho}_n(r(f), u)$$

*also satisfies*

$$P\ell_{\hat{f}} \leq \inf_{r \geq 1}(\mathcal{A}(r) + \kappa_2 \tilde{\rho}_n(r, u)),$$

*where*

$$\tilde{\rho}_n(r, u) = \rho_n\left(2r, u + \ln\frac{\pi^2}{6} + 2\ln(1 + c'_Y n + \log r)\right).$$

PROOF. The corollary follows directly from Theorems 4.11 and 2.2. We are able to remove the $\Theta^{2/p}$ term from Theorem 4.11 because $\Theta^{2/p} \leq c_{p,Y,A,\lambda}\frac{r^2}{n}$. $\quad\square$

The feature of this new bound that makes it better than our previous one is the fact that the term with the worst asymptotic behavior in $n$ has the best asymptotic behavior in $r$. Indeed, the $r^2$ term in $\rho_n(r, u)$ has a dependence on $n$ that scales like $1/n$, a much better rate than in the previous section. The significance of this is the suggestion that a regularization term of $\|f\|_H^2$ will result in over-regularization when $n$ is large. In fact, a study similar to the one at the end of Section 3 shows that Corollary 4.12 is indeed far better (we delay the details of this comparison until after Corollary 5.5, in which we improve the bound even further). In the following section, we will show that one can improve Corollary 4.12 even further by completely removing the $r^2$ term.

**5. Removing the $r^2$ term.**    The function $\rho_n$ from Corollary 4.12 is almost the function we would have liked to have. Its leading term is $\Theta^{2/(1+p)} \sim (r^{2p}n^{-1}\log^2 n)^{1/(1+p)}$, while the other term scales like $r^2/n$ and is dominant only for very large values of $r$. Here, we will show that the latter does not influence the minimization problem we are interested in and can be removed. Since some of the technical details of the proof of that observation are rather tedious and have already been presented in previous sections, certain parts of the argument will only be outlined.

Let us return to Theorem 2.2. The isomorphic condition we have established there holds in the set $F = rB_H$ with the functional

$$\psi(f, u) = c_{p,Y}\left(\max\{\Theta^{2/(1+p)}, \Theta^{2/p}\} + c_Y(1+u)\frac{\|f\|_\infty^2}{n}\right).$$

That is, for every $u > 0$, with probability at least $1 - \exp(-u)$, for every $f \in F$,

$$\tfrac{1}{2}P_n\mathcal{L}_f - \psi(f, u) \le P\mathcal{L}_f \le 2P_n\mathcal{L}_f + \psi(f, u).$$

Consider the minimization problem one faces when performing regularized learning. The problem is always to minimize a functional $\hat\Lambda = P_n\ell_f + \kappa_1 V_n$, hoping that the minimizer $\hat f$ will satisfy

$$P\ell_{\hat f} \le \inf_f \Lambda(f) = \inf_f(P\ell_f + \kappa_2 V_n),$$

where the functional $V_n : H \times \mathbb{R}_+ \to \mathbb{R}_+$ is nonnegative. In addition, all of the functionals we are interested in have the property that, for a fixed $f \in H$ and $u \in \mathbb{R}_+$, $V_n(f, u)$ tends to zero as $n \to \infty$.

We will specify our choice for the functional $V_n$ later, but, as a starting point, observe that since $f = 0$ is a potential minimizer, it follows that (assuming $\|Y\|_\infty \le 1$) any minimizer of $\hat\Lambda$ will satisfy $\hat\Lambda(\hat f) \le \hat\Lambda(0) \le 1 + V_n(0)$, and the same will hold for $\Lambda$. Since $V_n(0)$ tends to zero as $n$ grows, we can take $n$ sufficiently large (depending on $\|Y\|_\infty$) to ensure that $V_n(0) \le 1$. Therefore, for these values of $n$, any minimizer $\hat f$ of $\hat\Lambda$ satisfies

$$\hat\Lambda(\hat f) \le 2$$

and any minimizer $f^*$ of $\Lambda$ satisfies

$$\Lambda(f^*) \le 2.$$

Thus,

$$\{f : f \text{ minimizes } \Lambda\} \subset \{f : \mathbb{E}(f - Y)^2 \le 2\} \subset \{f : \mathbb{E}f^2 \le 9\}$$

and

$$\{f : f \text{ minimizes } \hat\Lambda\} \subset \{f : \hat\Lambda(f) \le 2\} \subset \{f : P_n f^2 \le 9\}.$$

Having this in mind, we will decompose $H$ into two subsets. The first, $H_1$, will contain $\{f : \mathbb{E}f^2 \leq 9\}$. In addition, we will show that $\bar{F}_r = H_1 \cap (r-1)B_H$ is an ordered, parameterized hierarchy of $H_1$ and that the assumptions of Theorem 2.5 will be satisfied with respect to a functional $V(r, x)$ for which the dominant term is $\Theta^{2/(1+p)}$.

Thus, by Theorem 2.5, with high probability, any minimizer of $\hat{\Lambda}$ in $H_1$ will satisfy

$$(5.1) \qquad P\ell_{\hat{f}} \leq \inf_{f \in H_1} (P\ell_f + \kappa_2 \tilde{V}(\|f\|_H, u)),$$

where $\tilde{V}$ is defined in a similar way to $\tilde{\rho}_n$ in Corollary 4.12.

The next step will be to extend the result beyond $H_1$ to $H$. Indeed, since $\{f : \mathbb{E}f^2 \leq 9\} \subset H_1$, the infimum in $H$ of the right-hand side of (5.1) is actually attained in $H_1$. Hence, the infimum in (5.1) is really over all functions in $H$. To conclude this line of reasoning, we will then show that with high probability, every empirical minimizer of $\hat{\Lambda}$ is in $H_1$, by proving that if $f \in H \setminus H_1$, then $P_n f^2 \geq 9$.

The correct decomposition of $H$ is attained using the following estimate for the ratio between the $\|f\|_H$ and $\|f\|_\infty$ for any function in $H$.

LEMMA 5.1. *Suppose that Assumptions 3.1 and 4.1 are satisfied. There is then a constant $\kappa_3 = \kappa_3(A, p, \|(\lambda_i)\|_{p,\infty})$ such that, for every $f \in H$,*

$$\mathbb{E}f^2 \geq \kappa_3 \left( \frac{\|f\|_\infty}{\|f\|_H^p} \right)^{2/(1-p)}.$$

PROOF. Recall that $\|K(x, x)\|_\infty \leq 1$ and let $r > 0$. Set $f(x) = \sum_{i=1}^\infty t_i \sqrt{\lambda_i} \times \varphi_i(x)$, where $\|t\|_2 = r$, and observe that $\|f\|_\infty \leq \|K\|_\infty r \leq r$. Also, since $\|(\lambda_i)\|_{p,\infty} < \infty$ and $(\lambda_i)_{i=1}^\infty$ is nonnegative and nonincreasing, it follows that for every $i$, $\lambda_i \leq (\|(\lambda_i)\|_{p,\infty}/i)^{1/p}$.

Fix $N$ (to be specified later) and observe that

$$\|f\|_\infty \leq A \left( \sum_{i=1}^N |t_i| \sqrt{\lambda_i} + r \left( \sum_{N+1}^\infty \lambda_i \right)^{1/2} \right)$$

$$\leq A \left( \sum_{i=1}^N |t_i| \sqrt{\lambda_i} + r \|(\lambda_i)\|_{p,\infty}^{1/2p} \left( \frac{1}{N} \right)^{(1-p)/2p} \right)$$

$$\leq A \sum_{i=1}^N |t_i| \sqrt{\lambda_i} + \frac{\|f\|_\infty}{2},$$

provided that $N^{(1-p)/2p} \geq 2Ar\|(\lambda_i)\|_{p,\infty}^{1/2p}/\|f\|_\infty$. Hence, $A \sum_{i=1}^N t_i \sqrt{\lambda_i} \geq \|f\|_\infty/2$. Note that $r/\|f\|_\infty$ is bounded below by 1 and so we can choose an

integer $N$ such that

$$\frac{2Ar\|(\lambda_i)\|_{p,\infty}^{1/2p}}{\|f\|_\infty} \leq N^{(1-p)/2p} \leq \frac{cAr\|(\lambda_i)\|_{p,\infty}^{1/2p}}{\|f\|_\infty}$$

for some constant $c$ depending on $p$ and $\|(\lambda_i)\|_{p,\infty}$. Clearly, for any $v \in \mathbb{R}^N$, $\|v\|_{\ell_2^N} \geq \|v\|_{\ell_1^N}/\sqrt{N}$ and thus,

$$\sum_{i=1}^N t_i^2 \lambda_i \geq c'\frac{\|f\|_\infty^2}{N} = c_1\left(\frac{\|f\|_\infty}{r^p}\right)^{2/(1-p)},$$

where $c_1$ is a constant depending on $A$, $p$ and $\|(\lambda_i)\|_{p,\infty}$.

On the other hand, since $(\varphi_i)_{i=1}^\infty$ is an orthonormal family, we have

$$\mathbb{E}f^2 = \mathbb{E}\sum_{i,j} t_i t_j \sqrt{\lambda_i \lambda_j}\varphi_i\varphi_j \geq \sum_{i=1}^N t_i^2 \lambda_i \geq c_1\left(\frac{\|f\|_\infty}{r^p}\right)^{2/(1-p)}. \qquad \square$$

Let

$$H_1 = \{0\} \cup \left\{f : \kappa_3\left(\frac{\|f\|_\infty}{\|f\|_H^p}\right)^{2/(1-p)} \leq 50\right\}.$$

Since the set of minimizers of any functional $\Lambda$ we will be interested in is contained in $\{f : \mathbb{E}f^2 \leq 9\}$, it follows, by Lemma 5.1, that the set of such minimizers is contained in $H_1$.

The set $H_1$ has additional properties. There is a constant $c$, depending on $p$ and $\kappa_3$, such that on $H_1$,

$$(5.2) \qquad\qquad\qquad \|f\|_\infty \leq c\|f\|_H^p.$$

Moreover, for every $r \geq 1$, if one considers $\bar{F}_r = H_1 \cap (r-1)B_H$, then the minimizer of $P\ell_f$ in $F_r = (r-1)B_H$ actually belongs to $\bar{F}_r$ (again, by comparing to $f = 0$). Therefore, it is straightforward to show that $\bar{F}_r$ is an ordered, parameterized hierarchy of $H_1$ with $r(f) = \|f\|_H + 1$, implying that one can obtain the desired isomorphic result on $H_1$, with the $\|f\|_\infty^2/n$ term replaced by $\|f\|_H^{2p}/n$.

Indeed, we can combine Theorem 2.2 with (5.2) and the fact that the localized averages $\mathbb{E}\|P_n - P\|$ indexed by $\{\text{star}(\mathcal{L}_{\bar{F}_r}, 0) : \mathbb{E}h \leq x\}$ are smaller than the localized averages indexed by the larger set $\{\text{star}(\mathcal{L}_{F_r}, 0) : \mathbb{E}h \leq x\}$ to show that for every $r \geq 1$, with probability at least $1 - \exp(-u)$, for every $f \in \bar{F}_r$,

$$\frac{1}{2}P_n\mathcal{L}_{r,f} - \frac{x}{2} - c(1+r^{2p})\frac{u}{n} \leq P\mathcal{L}_{r,f} \leq 2P_n\mathcal{L}_{r,f} + \frac{x}{2} + c(1+r^{2p})\frac{u}{n},$$

where $\mathcal{L}_{r,f}$ is the excess loss associated with $f$ relative to $\bar{F}_r$.

Using Theorem 4.11, one obtains the following result.

COROLLARY 5.2. *There exists a constant $\kappa_4'$ that depends on $p, A, \|(\lambda_i)\|_{p,\infty}$ and $\|Y\|_\infty$, for which the following holds. If $\Upsilon = r^p/\sqrt{n}$, then the function*

$$V'(r,u) = \kappa_4'(1+u)\max\{(\Upsilon\log n)^{2/(1+p)}, (\Upsilon\log n)^{2/p}, \Upsilon^2\}$$

*satisfies the hypothesis of Theorem 2.5 for the hierarchy $\{\bar{F}_r : r \geq 1\}$.*

*In particular, if we set $\hat{\Lambda}'(f,x) = P_n\ell_f + \kappa_1\tilde{V}'(f,u)$, then, with probability at least $1 - \exp(-u)$, every $f$ that minimizes $\hat{\Lambda}'$ in $H_1$ also satisfies*

$$P\ell_{\hat{f}} \leq \inf_{f \in H}(P\ell_f + \kappa_2\tilde{V}'(r(f),u)),$$

*where $\tilde{V}'$ is defined analogously to $\tilde{\rho}_n$ in Corollary 4.12.*

Next, we will show that the $(\Upsilon\log n)^{2/p}$ and $\Upsilon^2$ terms are nonessential. Indeed, for sufficiently large $n$, the minimal value in $H$ of $\hat{\Lambda}$ will be at most 2 (by comparing it to $f = 0$). Hence, if $f \in H$ satisfies $\kappa_5'\kappa_1\Upsilon\log n \geq 2$ [i.e., if $\|f\|_H \geq \kappa_5(n/\log^2 n)^{1/2p}$], then it is not a potential minimizer of $\hat{\Lambda}'$ in $H$. (Note that we can, by increasing $\kappa_4'$, take $\kappa_5$ as small as we like; this will be used later.) Therefore, on the set of potential minimizers, $\Upsilon\log n \leq c$, where $c$ depends on $\kappa_1$, $\kappa_4'$ and $p$. Hence, on this set of minimizers, we can bound

$$V'(r,u) \leq \kappa_4(1+u)(\Upsilon\log n)^{2/(1+p)}.$$

Denoting the right-hand side by $V(r,u)$, we can invoke Remark 2.6 to show that $V(r,u)$ is a valid functional.

Note that we can increase $H$ by adding every function $f \in H$ for which $\|f\|_H \geq (n/\log^2 n)^{(1/2p)}$; we have already argued that such functions cannot minimize $\hat{\Lambda}$.

To conclude, if

$$H_1' = H_1 \cup \{f : \|f\|_H \geq \kappa_5(n/\log^2 n)^{1/2p}\},$$

then, with probability at least $1 - \exp(-u)$, every $f$ that minimizes

$$P_n\ell_f + \kappa_1\tilde{V}(r(f),u)$$

in $H_1'$ also satisfies

$$P\ell_{\hat{f}} \leq \inf_{f \in H}(P\ell_f + \kappa_2\tilde{V}(r(f),u)).$$

Next, let us consider the set $H_2 = H \setminus H_1'$. Clearly, each function in $H_2$ satisfies $\|f\|_H \leq c_1\|f\|_\infty^{1/p}$ and $\mathbb{E}f^2 \geq 50$. We will show that, with high probability, any $f \in H_2$ satisfies $P_nf^2 \geq 9$ and thus is not a potential minimizer of $\hat{\Lambda}$ in $H$.

Lemma 5.3. *There exists a constant $\kappa_6$ that depends on $A$, $p$, and $\|(\lambda_i)\|_{p,\infty}$ and an absolute constant $\kappa_7$ for which the following holds. If $0 \in F$ and $F \subset \kappa_6(n/\log^2 n)^{1/2p}B_H$, then, for every $u > 0$, with probability at least $1 - \exp(-u)$, for every $f \in F$,*

$$P_n f^2 \geq \frac{1}{2}\mathbb{E}f^2 - 1 - \kappa_7(1 + \|F\|_\infty^2)\frac{u}{n},$$

*where $\|F\|_\infty = \sup_{f \in F}\|f\|_\infty$.*

Proof.   Apply Theorem 4.11 with $Y \equiv 0$, noting that, in this case, $\mathcal{L}_f = f^2$. It follows that we can set

$$W_{x,r} = \{f^2 : \|f\|_H \leq r, \mathbb{E}f^2 \leq x\}$$

and $\mathbb{E}\|P_n - P\|_{W_{x,r}} \leq x/8$, provided that

$$x \geq c_1 \max\{(\Upsilon \log n)^{2/(1+p)}, (\Upsilon \log n)^{2/p}\},$$

where $c_1$ depends on $A$, $p$ and $\|(\lambda_i)\|_{p,\infty}$. We will apply this fact for $x = 2$. That is, we need to ensure that $r$ is chosen in such a way that

$$c_1 \max\{(\Upsilon \log n)^{2/(1+p)}, (\Upsilon \log n)^{2/p}\} \leq 2,$$

which is the case, for example, if $r \leq c_2(n/\log^2 n)^{1/2p}$.

The result now follows from Theorem 2.2.   □

Set $r_H = \kappa_6(n/\log^2 n)^{1/2p}$ and recall that $\kappa_5$ can be taken as small as we like. In particular, we may assume that $\kappa_5 \leq \kappa_6$ and so $H_2 \subset r_H B_H$.

The final preliminary step we take is to decompose $H_2$ into $L_\infty$-shells in the following way. Fix $u > 0$ and set $r_0$ such that $\kappa_7 u(1 + r_0^2)/n < 9$. Define $(r_i)_{i=0}^m$ by $r_i = 2^i r_0$, where $m$ is the smallest number such that $r_m \geq r_H$. Thus, $m \leq c_1(\log n + \log u)$. Let

(5.3)        $$B = \left\{f : \|f\|_\infty \geq \kappa_8\|f\|_H\left(\frac{u}{n}\right)^{(1-p)/2p}\right\},$$

where $\kappa_8$ is some constant that will be named in the proof of the following lemma. We will consider the sets $F_0 = H_2 \cap r_0 B_\infty$ and

$$F_i = H_2 \cap \{f : r_i \leq \|f\|_\infty \leq r_{i+1}\} \cap B.$$

Since $\bigcup_{i=0}^m (H_2 \cap \{f : r_i \leq \|f\|_\infty \leq r_{i+1}\}) = H_2$, any $f \in H_2 \setminus \bigcup_{i=0}^m F_i$ satisfies

$$\|f\|_\infty \leq \kappa_8\|f\|_H\left(\frac{u}{n}\right)^{(1-p)/2p}$$

and because $\|f\|_H \leq r_H$, we have

$$\|f\|_\infty \leq \kappa_6\kappa_8\left(\frac{n}{\log^2 n}\right)^{1/2p} \cdot \left(\frac{u}{n}\right)^{(1-p)/2p} = c_1 u^{(1-p)/2p}\frac{n^{1/2}}{\log^{1/p} n}.$$

Therefore,

$$\frac{\|H_2 \setminus \bigcup_{i=0}^m F_i\|_\infty^2}{n} \leq c_1^2 \frac{u^{(1-p)/p}}{\log^{2/p} n}.$$

LEMMA 5.4. *There exist constants $c_1$ and $c_2$, depending only on $A$, $p$ and $\|(\lambda_i)\|_{p,\infty}$, for which the following holds. Fix $n$ and $0 < u < c_1 n$, and perform the above decomposition. For every $0 \leq i \leq m$, with probability at least $1 - \exp(-u)$, every $f \in F_i$ satisfies $P_n f^2 \geq 9$. Also, if $u \leq c_2 (\log n)^{2/(1-p)}$, then, with probability $1 - \exp(-u)$, for every $f \in H_2 \setminus \bigcup_{i=0}^m F_i$, $P_n f^2 \geq 9$.*

PROOF. First, fix $1 \leq i \leq m$ and apply Lemma 5.3 to the set $F_i$. For every $f \in F_i$, $\|f\|_\infty \leq \|F_i\|_\infty \leq 2\|f\|_\infty$, and thus, with probability at least $1 - \exp(-u)$,

$$P_n f^2 \geq \frac{1}{2}\mathbb{E}f^2 - 1 - \kappa_7 \frac{u(1 + \|F_i\|_\infty^2)}{n} \geq \frac{1}{2}\mathbb{E}f^2 - 1 - 2\kappa_7 \frac{u(1 + \|f\|_\infty^2)}{n}.$$

On the other hand, for every $f \in B$,

$$\frac{1}{4}\mathbb{E}f^2 \geq \frac{\kappa_3}{4}\left(\frac{\|f\|_\infty}{\|f\|_H^p}\right)^{2/(1-p)} \geq 2\kappa_7 \frac{u\|f\|_\infty^2}{n},$$

provided that $\kappa_8 \geq (8\kappa_7/\kappa_3)^{(1-p)/2p}$. Therefore, with probability at least $1 - \exp(-u)$, for every $f \in F_i$,

$$P_n f^2 \geq \frac{1}{4}\mathbb{E}f^2 - 1 - \frac{2\kappa_7 u}{n} \geq 10 - \frac{2\kappa_7}{c_1} \geq 9$$

for a suitably large choice of $c_1$.

Turning to $F_0$, since $\kappa_7 \frac{u(1+\|F_0\|_\infty^2)}{n} \leq 9$, we have, by Lemma 5.3, with probability at least $1 - \exp(-u)$, for every $f \in F_0$,

$$P_n f^2 \geq \frac{1}{2}\mathbb{E}f^2 - 1 - \kappa_7 \frac{u(1+r_0^2)}{n} \geq 9.$$

Finally, since $n^{-1}\|H_2 \setminus \bigcup_{i=0}^m F_i\|_\infty^2 \leq c\frac{u^{(1-p)/p}}{\log^{2/p} n}$, it follows that for our choice of $u$,

$$\kappa_7 u \frac{\|H_2 \setminus \bigcup_{i=0}^m F_i\|_\infty^2}{n} \leq 9$$

from which our claim follows, using the same argument as for $F_0$. □

We can now prove our main result, which is the second part of the following claim and was formulated as Theorem A in the Introduction.

COROLLARY 5.5. *If Assumptions 3.1 and 4.1 are satisfied, then there exist constants $c_1$, $c_2$ and $c_3$ that depend only on $A$, $p$ and $\|(\lambda_i)\|_{p,\infty}$, a constant $N_0$ that depends on $\|Y\|_\infty$ and $p$ and a constant $c_Y$ that depends only on $\|Y\|_\infty$, for which the following holds.*

*If $n \geq N_0$, $c_1 \log \log n \leq u \leq c_2 (\log n)^{2/(1-p)}$, then, with probability at least $1 - \exp(-u/2)$, for every $f \in H_2$, $P_n f^2 \geq 9$. Thus, all of the minimizers in $H$ of*

$$(5.4) \qquad P_n \ell_f + \kappa_1 \tilde{V}(f, u)$$

*belong to $H_1$. In particular, for such values of $u$, with probability at least $1 - 2\exp(-u/2)$, every minimizer $\hat{f}$ in $H$ of (5.4) satisfies*

$$P\ell_{\hat{f}} \leq \inf_{f \in H} (P\ell_f + \kappa_2 \tilde{V}(f, u)),$$

*where*

$$\tilde{V}(f, u) = c_3 (1 + u + c_Y \ln n + \ln \log(\|f\|_H + e)) \left( \frac{(\|f\|_H + 1)^p \log n}{\sqrt{n}} \right)^{2/(1+p)}.$$

Let us (briefly) repeat the analysis that we carried out at the end of Section 3. Recall Assumption 3.2: we assume the existence of $0 < \sigma < 1$ such that the regression function $\mathbb{E}(Y|X)$ belongs to $T_K^\sigma L_2$. Recall, also, that under this assumption, the approximation error $\mathcal{A}(r)$ behaves like $r^{-4\sigma/(1-2\sigma)}$. Under Corollary 5.5, the error of the empirical minimizer is like (ignoring logarithmic terms)

$$\inf_{r \geq 1} \left( \mathcal{A}(r) + \frac{r^{2p/(1+p)}}{n^{1/(1+p)}} \right) \lesssim \inf_{r \geq 1} \left( \frac{1}{r^{4\sigma/(1-2\sigma)}} + \frac{r^{2p/(1+p)}}{n^{1/(1+p)}} \right),$$

which can be optimized by choosing $r = n^{-k/(2p+kp+k)}$. This gives us a final error rate of

$$\left( \frac{1}{n} \right)^{2\sigma/(p+2\sigma)},$$

which is, as promised, better by a polynomial factor than the previous error rate of $n^{-2\sigma/(p+1)}$ whenever $\sigma < 1/2$.

## APPENDIX: PROOFS

The starting point in the proof of Theorem 2.5 is the following theorem by Bartlett [1].

THEOREM A.1. *Suppose that $\{F_r; r \geq 1\}$ is an ordered, parameterized hierarchy and that $\rho_n(r)$ is a positive, continuous, increasing function. If, for all $r \geq 1$ and all $f \in F_r$,*

$$(A.1) \qquad \tfrac{1}{2} P_n \mathcal{L}_{r,f} - \rho_n(r) \leq P\mathcal{L}_{r,f} \leq 2P_n \mathcal{L}_{r,f} + \rho_n(r),$$

*then*

$$P\ell_{\hat{f}} \leq \inf_{f \in F}(P\ell_f + c_1\rho_n(r(f))),$$

*where $\hat{f}$ is any function that minimizes the functional $P_n\ell_f + c_2\rho_n(r(f))$.*

PROOF. Let $(r_i)_{i=1}^{\infty}$ be an increasing sequence (to be determined later) such that $r_1 = 1$ and $r_i \to \infty$ as $i \to \infty$. Define, for each $i \geq 1$, $u_i = u + \ln(\pi^2/6) + 2\ln i$. Then

$$\sum_{i=0}^{\infty} e^{-u_i} = e^{-u}$$

and so, by the union bound, with probability at least $1 - e^{-u}$, for every $i \geq 1$,

$$\tfrac{1}{2}P_n\mathcal{L}_{r_i,f} - \rho_n(r_i, u_i) \leq P\mathcal{L}_{r_i,f} \leq 2P_n\mathcal{L}_{r_i,f} + \rho_n(r_i, u_j).$$

If we only cared about a sequence of $r_i$, this would be enough for our result. However, we need an almost-isomorphic condition for all $r \geq 1$ and so the next step must be to find an almost-isomorphic condition for $F_r$ when $r \in [r_{j-1}, r_j]$. In one direction, we have

$$\begin{aligned}
P\mathcal{L}_{r,f} &= P\mathcal{L}_{r_j,f} - P\mathcal{L}_{r_j,f_r^*} \\
&\leq 2P_n\mathcal{L}_{r_j,f} + \rho_n(r_j, u_j) - P\mathcal{L}_{r_j,f_r^*} \\
\text{(A.2)} \quad &= 2P_n\mathcal{L}_{r,f} + 2P_n\mathcal{L}_{r_j,f_r^*} + \rho_n(r_j, u_j) - P\mathcal{L}_{r_j,f_r^*} \\
&\leq 2P_n\mathcal{L}_{r,f} + 5\rho_n(r_j, u_j) + 3P\mathcal{L}_{r_j,f_r^*} \\
&\leq 2P_n\mathcal{L}_{r,f} + 5\rho_n(r_j, u_j) + 3P\mathcal{L}_{r_j,f_{r_{j-1}}^*},
\end{aligned}$$

while in the other direction, we get

$$\begin{aligned}
2P\mathcal{L}_{r,f} &= 2P\mathcal{L}_{r_j,f} - 2P\mathcal{L}_{r_j,f_r^*} \\
&\geq P_n\mathcal{L}_{r_j,f} - 2\rho_n(r_j, u_j) - 2P\mathcal{L}_{r_j,f_r^*} \\
\text{(A.3)} \quad &= P_n\mathcal{L}_{r,f} + P_n\mathcal{L}_{r_j,f_r^*} - 2\rho_n(r_j, u_j) - 2P\mathcal{L}_{r_j,f_r^*} \\
&\geq P_n\mathcal{L}_{r,f} - \tfrac{5}{2}\rho_n(r_j, u_j) - \tfrac{3}{2}P\mathcal{L}_{r_j,f_r^*} \\
&\geq P_n\mathcal{L}_{r,f} - \tfrac{5}{2}\rho_n(r_j, u_j) - \tfrac{3}{2}P\mathcal{L}_{r_j,f_{r_{j-1}}^*}.
\end{aligned}$$

We can now choose our sequence $r_i$: recall that $r_1 = 1$ and set $r_i$, for all $i \geq 2$, to be the largest number satisfying both of the following inequalities:

$$\begin{aligned}
\text{(A.4)} \quad & r_i \leq 2r_{i-1}, \\
& P\mathcal{L}_{r_i,f_{r_{i-1}}^*} \leq \rho_n(r_i, u_i).
\end{aligned}$$

Note that choosing the largest number is not a problem because both $\rho_n(r, u)$ and $P\mathcal{L}_{r, f^*_{r_{j-1}}}$ are continuous functions of $r$; that is, the supremum of the set of $r$ satisfying (A.4) is attained.

Our choice of $r_i$ ensures that, for all $i \geq 1$,

$$(A.5) \qquad i \leq \frac{P\ell(f^*_{r_1}, Y)}{\rho_n(r_1, u_1)} - \frac{P\ell(f^*_{r_i}, Y)}{\rho_n(r_i, u_i)} + \log(2r_i) \leq \frac{P\ell(f^*_{r_1}, Y)}{\rho_n(r_1, u_1)} + \log(2r_i).$$

Indeed, for $i = 1$, this is trivial. For larger $i$, we can proceed by induction: our definition of $r_i$ ensures that either $r_i = 2r_{i-1}$ or $P\ell(f^*_{r_{i-1}}, Y) = P\ell(f^*_{r_i}, Y) + \rho_n(r_i, u_i)$. In the first case, $\log r_i = \log r_{i-1} + 1$ and the inductive step follows. In the second case, assuming that

$$i - 1 \leq \frac{P\ell(f^*_{r_1}, Y)}{\rho_n(r_1, u_1)} - \frac{P\ell(f^*_{r_{i-1}}, Y)}{\rho_n(r_{i-1}, u_{i-1})} + \log r_{i-1},$$

it follows that

$$\begin{aligned}
i &\leq \frac{P\ell(f^*_{r_1}, Y)}{\rho_n(r_1, u_1)} - \frac{P\ell(f^*_{r_{i-1}}, Y)}{\rho_n(r_{i-1}, u_{i-1})} + 1 + \log(2r_i) \\
&\leq \frac{P\ell(f^*_{r_1}, Y)}{\rho_n(r_1, u_1)} - \frac{P\ell(f^*_{r_{i-1}}, Y)}{\rho_n(r_i, u_i)} + 1 + \log(2r_i) \\
&= \frac{P\ell(f^*_{r_1}, Y)}{\rho_n(r_1, u_1)} - \frac{P\ell(f^*_{r_i}, Y)}{\rho_n(r_i, u_i)} + \log(2r_i),
\end{aligned}$$

which proves (A.5) by induction. In particular, for any $i \geq 1$ and any $r \geq r_i$, $u_i \leq \theta(r, u)$. Therefore,

$$\rho_n(r_i, u_i) \leq \rho_n(2r, \theta(r, u))$$

for any $r \in [r_{i-1}, r_i]$.

Note that (A.5) implies that the sequence $r_i$ tends to infinity with $i$. Then, by (A.2), (A.3) and (A.4), with probability at least $1 - e^{-u}$, for all $r \geq 1$ and all $f \in F_r$,

$$\tfrac{1}{2} P_n \mathcal{L}_{r, f} - 4\rho_n(2r, \theta(r, u)) \leq P\mathcal{L}_{r, f} \leq 2P_n \mathcal{L}_{r, f} + 8\rho_n(2r, \theta(r, u)).$$

We conclude the proof by applying Theorem A.1.  □

CENTRE FOR MATHEMATICS
  AND ITS APPLICATIONS
THE AUSTRALIAN NATIONAL UNIVERSITY
CANBERRA, ACT 0200
AUSTRALIA
AND
DEPARTMENT OF MATHEMATICS
TECHNION, I.I.T.
HAIFA, 32000
ISRAEL
E-MAIL: shahar.mendelson@anu.edu.au

DEPARTMENT OF STATISTICS
UNIVERSITY OF CALIFORNIA
BERKELEY, CALIFORNIA 94720
USA
E-MAIL: joeneeman@gmail.com